\documentclass[reqno, 12pt]{amsart}
\def\C{{\mathbb C}}

\def\Q{{\mathbb Q}}

\def\Z{{\mathbb Z}}
\def\Zp{{\mathbb Z}_p}
\def\R{{\mathbb R}}

\def\Fp{{{\mathbb F}_p}}
\def\F{{\mathbb F}}

\def\Fr{{{\mathbb F}_r}}

\def\dis{\displaystyle}

\theoremstyle{definition}
\newtheorem{defn}{Definition}   

\newtheorem{question}{Question}

\theoremstyle{remark}
\newtheorem{rem}{Remark}        
      
\newtheorem{example}{Example}

\begin{document}
\title[Can a Drinfeld module be modular?]
{Can a Drinfeld module be modular?}
\author{David Goss}
\address{Department of Mathematics\\ The Ohio State University\\ 231 W.
$18^{\text{th}}$ Ave. \\ Columbus, Ohio 43210}
\email{goss@math.ohio-state.edu}
\date{December 23, 2002}

\begin{abstract}
Let $k$ be a global function field with field of constants $\Fr$,
$r=p^m$, and let
$\infty$ be a fixed place of $k$.
In his habilitation thesis \cite{boc2}, Gebhard B\"ockle attaches 
abelian Galois representations to characteristic $p$ valued
cusp eigenforms and double cusp eigenforms \cite{go1} 
such that Hecke eigenvalues
correspond to the image of Frobenius elements.  In the case where
$k=\Fr(T)$ and $\infty$ corresponds to the pole of $T$, it then becomes
reasonable to ask whether rank $1$ Drinfeld modules over $k$ are themselves
``modular'' in that their Galois representations arise from a cusp or
double cusp form. This paper gives an introduction to \cite{boc2}
with an emphasis on modularity and closes with some specific questions
raised by B\"ockle's work.
\end{abstract}

\maketitle

\section{Introduction}\label{intro}
Let $k$ be a number field and let $E$ be an ``arithmetic object'' over
$k$ such as an elliptic curve or abelian variety. Following Riemann, 
Artin, Weil, Hasse, and Grothendieck, one associates to $E$ an $L$-series
$L(E,s)$ via its associated Galois representations. Thus for
each finite prime $\mathfrak p$ of $k$, one obtains (or is conjectured to
obtain) a canonical polynomial $f_{\mathfrak p}(u)\in \Z[u]$ and 
one sets
$$L(E,s)=\prod_{\mathfrak p}f_{\mathfrak p}(N{\mathfrak p}^{-s})^{-1}\,.$$
Using
estimates, such as those arising from the Weil Conjectures, one sees that
this Euler product converges on a non-trivial half-plane of the
complex numbers $\C$ to an analytic function. 

Uncovering the properties of $L(E,s)$ is a major goal of number theory.
The primary approach to this end, also initiated by Riemann, is to
equate $L(E,s)$ with ``known''  or ``standard''
Dirichlet series via a reciprocity law.
For instance, as recalled in Subsection \ref{thetadir}, 
the Riemann zeta function completed with $\Gamma$-factors at the
infinite primes, can also be obtained via an integral transform
of a theta-function; the analytic properties of the zeta function are then
consequences of those of the theta function. In general, for arbitrary $E$,
one may, conjecturally(!!) work the same way by replacing 
the theta-function with an ``automorphic form.'' In this fashion, one hopes
to show that the function $L(E,s)$ has an analytic continuation and
a functional equation under $s\mapsto k-s$ for some integer $k$.

The profundity of the task of attaching an appropriate automorphic
form to $L(E,s)$ is attested to by noting that Fermat's Last Theorem
follows as a consequence when $E$ is restricted to just the set of
semi-stable elliptic curves over $\Q$ \cite{wi1}. 

Now let $k$ be a global function field over a finite field $\Fr$,
$r=p^m$. Beginning with E.\ Artin's thesis, number theorists learned
how to attach $L$-series to arithmetic objects over $k$.
Grothendieck \cite{gro1} presented a cohomological
approach to these $L$-series which showed that they possess
an analytic continuation (as a rational function in $u=r^{-s}$) and
a functional equation of classical type. Moreover, the notion of automorphic
form is supple enough to work over $k$ also. Thus a very natural,
and important, problem was to investigate  whether the class of
$L$-series associated $k$ would also be included in the class of standard
$L$-series. With the recent work of L.\ Lafforgue (\cite{laf1}, see also
\cite{lau1}) this has now been established.

Lafforgue's work builds on the ideas and constructions of
V.G.\ Drinfeld and, in particular, his notion of an ``elliptic
module'' \cite{dr1}. Elliptic modules (now called ``Drinfeld modules'')
are analogs of elliptic curves and abelian varieties. However, they
are not projective objects; rather they are linear objects equipped with an
exotic action by an affine sub-algebra of $k$. More precisely, following
Drinfeld one picks a place of $k$, labels it ``$\infty$,'' and then sets
$A$ to be those elements in $k$ which have no poles away from $\infty$.
The ring $A$ then becomes analogous to $\Z$ and the field $k$
analogous to $\Q$. What makes this analogy especially convincing is that
Drinfeld modules arise over the algebraic closure of the
completion $k_\infty$ via
``lattices'' and ``exponential functions'' in a fashion
rather analogous to what
happens over the complex numbers
with the classical exponential function and elliptic curves. 
Thus, for instance, the moduli spaces of
Drinfeld modules of a given rank have \emph{both} an algebraic
and analytic description. Using the moduli curves of Drinfeld modules
of rank $2$, corresponding to rank $2$ lattices, Drinfeld (ibid.)
established his first general $2$-dimensional reciprocity law. This then
implies that elliptic curves with split-multiplicative
reduction at $\infty$ are
isogenous to Jacobian factors of these curves 
(see Section \ref{ellmodffields}).

The analogy between Drinfeld modules and elliptic curves extends far beyond
just the way these objects are constructed. Indeed, like an elliptic
curve, one can associate to a Drinfeld module $E$ its Tate modules
with their canonical Galois representations and Frobenius actions.
Thus, as with elliptic curves, it is completely natural to encode
this information into a characteristic $p$ valued $L$-function 
$L(E,s)$ where $s$ lies in the space $S_\infty$ (see Equation (\ref{sinfty})). 
Moreover, as with elliptic curves, Drinfeld modules naturally have a
theory of modular forms (defined in almost exactly the classical manner)
associated to them \cite{go1}. In the case where the rank is $2$, these modular
forms naturally live on an ``algebraic upper half-plane'' which plays the
role of the usual complex upper half-plane. Furthermore, these modular
forms come equipped with an action of the ``Hecke operators''
which are again defined following classical theory. However,
the relationship between the Hecke eigenvalues and the ``$q$-expansion''
coefficients of a given eigenform was, and is, very mysterious unlike
classical theory where it is quite transparent. 

In classical theory the parameter $q$ at $\infty$ satisfies
$dq=*qdz$. In the characteristic $p$ theory it satisfies 
$dq=*q^2dz$. As such, one is led to study both cusp forms (forms which vanish
at all cusps) and the subspace of ``double-cusp'' forms (cusp forms which also
have first order vanishing); both of these are readily seen to be
Hecke-modules.

In this paper we report on the seminal work \cite{boc2} of
Gebhard B\"ockle in which Galois representations are naturally associated
to cusp and double cusp forms. Previous to this work, B\"ockle and
R.\ Pink had developed a good cohomology theory associated to
``$\tau$-sheaves'' (which are a massive generalization of Drinfeld modules and
form the correct category in which to discuss characteristic $p$
valued $L$-series). Indeed, in \cite{boc1} B\"ockle used this cohomology
to establish very generally good analytic properties for these $L$-functions
of $\tau$-sheaves.
In \cite{boc2} the author applies the full power of this cohomological
theory to cusp forms associated to rank $2$
Drinfeld modules via the $\tau$-sheaves naturally associated to the
universal families lying over the moduli curves. By comparison with
similar constructions in the \'etale topology, the associated Galois
representations emerge. As the Hecke operators $T(I)$ in characteristic $p$
satisfy $T(I^2)=T(I)^2$ for \emph{all} $I$, one sees that the simple
Galois factors are \emph{abelian}. For cusp forms which are \emph{not}
double cuspidal, these representations essentially arise from \emph{finite}
abelian extensions (split totally at $\infty$) but for most double
cusp forms the representations have infinite image.

Classical theory immediately 
leads to an immense number of interesting questions
about these representations and their associated $L$-series (which indeed
have good analytic 
properties via the techniques in \cite{boc1}). As of now, one
does not even have good guesses as to what the answers might be.

This paper is written in
order to motivate interested number theorists to become involved
in these basic issues. While \cite{boc2} is daunting in the
great number of details that must be checked, this paper will be
quite short on details. Rather we focus on the ``big picture'' of 
how the characteristic $p$ theory compares with classical theory for
both number fields and function fields. For ease of exposition we
let $k=\Fr(T)$ and $\infty$ the place associated the pole of $T$ as
usual. Because the class number of $k$ (in terms of divisors of degree
$0$) is $1$, there exist many Drinfeld modules of rank $1$ defined over
$k$ and it now makes sense to ask if any of them are ``modular'' in that
their Galois representations arise from cusp forms.

In fact, while we now know that a Hecke eigenform $f$ gives rise to a
good $L$-series $L(f,s)$, we have no idea yet how to classify the functions
which arise nor do we know any sort of ``converse'' theorems. Still, it makes
sense to broaden the definition of ``modularity'' in order to allow one
to capture the $L$-series of the rank $1$ Drinfeld module \emph{up to
translation} (much as $\zeta(2s)$ is naturally associated to the
classical theta function; see Subsection \ref{classmod} just below). We then
find that there are really \emph{two} distinct notions of modularity depending
on whether the cusp form is double cuspidal or not. We will see that
the Carlitz module then becomes modular in both senses. Finally,
in Subsection \ref{rank1taumod} we
present a certain rank one Drinfeld module $C^{(-\theta)}$ (with
$C^{(-\theta)}_T(x):=\theta x-\theta x^r$)
whose associated Galois
representations quite conceivably --- with our current knowledge --- might
arise directly from a double cusp form. Classical theory certainly implies
that the answer as to whether this Drinfeld module is truly modular
or not should be very interesting.

It is my great pleasure to thank Gebhard B\"ockle for his immense patience
in guiding me through his thesis. Without his careful answers to my
many questions this paper would have been impossible. Indeed,
it is my sincere hope
that this work makes \cite{boc2} more accessible. Still,
any mistakes in this paper are the fault of its author. This paper
is an expanded version of a lecture presented at the 
Canadian Number Theory Association in May, 2002. It is also my pleasure
to thank the Association for the opportunity to present these ideas.
I also thank A.\ Greenspoon, D.\ Rohrlich and J.-P.\ Serre for their
help with earlier versions of this work.
Finally, I am very grateful to M.\ Ram Murty for suggesting that
I write an exposition based on my presentation.

\section{Classical Modularity over $\Q$}\label{classmod}

\subsection{Theta functions and Dirichlet characters}\label{thetadir}
The connection between modular forms and $L$-series is a central
theme of modern number theory. We will summarize some of the relevant
ideas in this section. An excellent source in this regard is \cite{kn1} which
we follow rather closely.

The theory begins with
Riemann's original paper on the distribution of primes.
Indeed, let
\begin{equation}\label{theta0}
\theta(z):=\sum_{n=-\infty}^\infty e^{in^2\pi z}
\end{equation}
be the classical theta function. One knows that
$\theta (z)$ is analytic on the upper half-plane
${\mathcal H}:=\{z \in \C\mid \Im (z)>0\}$, and
one visibly sees that
\begin{equation}\label{theta1}
\theta(z+2)=\theta(z)\,.
\end{equation}
Due to Jacobi and Poisson, one also also has
the much deeper formula
\begin{equation} \label{theta2}
\theta(-1/z)=(z/i)^{1/2}\theta(z)\,,
\end{equation}
where, for the square root,
one takes the principal value which is cut on the negative real axis. 
Transformation laws (\ref{theta1}) and (\ref{theta2}) are summarized
by saying that $\theta(z)$ is a {\it modular form of weight $1/2$}
associated to the group $\Gamma_\theta$ of automorphisms of
$\mathcal H$ generated by $z\mapsto z+2$ and $z \mapsto -1/z$.

The well-known application of $\theta (z)$ to $L$-series, due to
Riemann, then arises in the following fashion.
Let 
\begin{equation}\label{zeta0}
\zeta(s):=\sum_{n=1}^\infty n^{-s}=\prod_{p~\rm prime}(1-p^{-s})^{-1}\,,
\end{equation}
be the Riemann zeta function
and let $\dis\Gamma (s)=\int_0^\infty t^se^{-t}\, \frac{dt}{t}$ be Euler's 
gamma function. One sets $\dis \Lambda(s):=\zeta(s)\Gamma(s/2)\pi^{-s/2}$.
Through the integral representation for $\Gamma(s)$ and
a change of variables one finds
\begin{equation}\label{zeta1}
2\Lambda(s)=\int_0^\infty \left(\theta(it)-1\right)
t^{s/2}\, \frac{dt}{t}\,.\end{equation}
Equation (\ref{zeta1}) is, in turn, rewritten as
\begin{equation}\label{zeta2}
2\Lambda(s)=\int_0^1\theta(it)t^{s/2}\,\frac{dt}{t}-
\int_0^1t^{s/2}\,\frac{dt}{t}+\int_1^\infty
\left(\theta(it)-1\right)t^{s/2}\, \frac{dt}{t}\,.
\end{equation}
The second term on the right is readily computed to equal $2/s$.
Via Equation (\ref{theta2}), the first term on the right
in Equation (\ref{zeta2}) is computed
to be $\dis \int_1^\infty (\theta(it)-1)t^{1/2(1-s)}\, \frac{dt}{t}-\frac{2}{1-s}\,.$ Thus, finally,
\begin{equation}\label{zeta3}
2\Lambda(s)=\int_1^\infty (\theta(it)-1)t^{s/2}\, \frac{dt}{t}+
\int_1^\infty (\theta(it)-1)t^{1/2(1-s)}\,\frac{dt}{t}-\frac{2}{s(1-s)}\,.
\end{equation}
The first two terms on the right can be shown to be entire in $s$. 
Moreover, from the invariance of the right hand side of (\ref{zeta3})
under $s\mapsto 1-s$,
we deduce that
\begin{equation}\label{zeta4}
\Lambda(s)=\Lambda(1-s)\,.
\end{equation}
This is the famous {\it functional equation} for $\zeta(s)$, and surely
one of the most sublime statements in mathematics.

\begin{rem}\label{remark1}
The above argument actually gives {\it both} the analytic continuation
of $\zeta(s)$ and the functional equation (\ref{zeta4}). In the following
we will use ``functional equation'' to mean both an analytic continuation
and invariance under $s\mapsto k-s$ for some integer $k$.
\end{rem}

\subsubsection{$L$-series associated to modular forms} \label{modularforms}

The derivation of the functional equation (\ref{zeta4}) of $\zeta(s)$ from the
properties (\ref{theta1}), (\ref{theta2}) of $\theta(z)$ is just
the very tip of the iceberg as we shall see. We begin by recalling
the general definition of a modular form.

Let $\gamma:=\dis \left(\begin{array}{cc}a&b\\c&d\end{array}\right)\in SL_2(\R)$ and let $z\in\mathcal H$. We set
\begin{equation}\label{gammaact}
\gamma z:=\frac{az+b}{cz+d}\,.\end{equation}
The map $z \mapsto \gamma z$ is clearly
an analytic automorphism of $\mathcal H$ (the inverse automorphism arising
from the inverse matrix). Note that $\gamma$ and $-\gamma$ have the same
action on $\mathcal H$. 
Let $\Gamma:=SL_2(\Z)$ and let $\tilde{\Gamma}$ be a subgroup of
$\Gamma$. Finally, let $k$ be a real number and assume that we have
chosen a branch so that $z^k$ is analytic on $\mathcal H$.

\begin{defn} \label{unrestmod}
Let $f(z)$ be an analytic function on $\mathcal H$. We say that
$f(z)$ is an {\it unrestricted modular form of weight $k$ associated to
$\tilde{\Gamma}$} if and only if
\begin{equation} \label{modular1}
\hbox{$f(\gamma z)=(cz+d)^kf(z)$\quad for \quad
$\left(\begin{array}{cc}a&b\\c&d\end{array}\right)\in\tilde{\Gamma}$.}
\end{equation}
\end{defn}

\noindent 
More generally, one allows certain constants, called ``multiplier systems,''
in front of $(cz+d)^k$ in Equation (\ref{modular1}). For instance, multipliers
are needed in order for $\theta(z)$ to be modular as in Equation
(\ref{theta2}). 

Now, let $\chi$ be a primitive
Dirichlet character modulo $m$. To $\chi$ one associates
the $L$-series
\begin{equation}\label{dirichlet1}
L(\chi,s):=\sum_{n=1}^\infty \chi(n)n^{-s}=\prod_{p~\rm prime}(1-\chi(p)
p^{-s})^{-1}\,,\end{equation}
where, by definition, $\chi(n)=0$ if and only if $\gcd (n,m)>1$.
(So, if $\chi=\chi_0$ the trivial character, then $L(\chi,s)=\zeta(s)$.)
To $\chi$ one also associates a theta function as follows. Let
$z\in \mathcal H$ and set
\begin{equation}\label{theta3}
\theta(\chi,z):=
\begin{cases}\sum_{n=-\infty}^\infty \chi(n)e^{in^2\pi z /m}
& \text{if $\chi(-1)=1$}\\
\sum_{n=-\infty}^\infty \chi(n) n e^{in^2\pi z /m}& \text{if $\chi(-1)=-1$\,.}
\end{cases}
\end{equation}

It is clear that $\theta(\chi_0,z)=\theta(z)$, and that
\begin{equation}\label{theta4}
\theta(\chi,z+2m)=\theta(\chi,z)\,.
\end{equation}
Moreover, 
\begin{equation}\label{theta5}
\theta(\chi,-1/z):=\begin{cases} w(\chi,m)(z/i)^{1/2}\theta(\bar{\chi},z)
&\text{if $\chi(-1)=1$}\\
-iw(\chi,m)(z/i)^{3/2}\theta(\bar{\chi},z)& \text{if $\chi(-1)=-1$\,,}
\end{cases}\end{equation}
where $|w(\chi,m)|=1$ and $\bar{\chi}$ is the complex conjugate
character. In particular, if $\chi=\bar{\chi}$, we obtain a modular
form (of weights $1/2$ or $3/2$) for the group $\Gamma_{\theta,\chi}$
of automorphisms
of $\mathcal H$ generated by $z\mapsto -1/z$ and
$z \mapsto z+2m$.

One now sets
\begin{equation}\label{lserieschi}
\Lambda(\chi,s):=\begin{cases} m^{s/2}\Gamma(s/2)\pi^{-s/2}L(\chi,s)
&\text{if $\chi(-1)=1$}\\
m^{(s+1)/2}\Gamma\left(\frac{s+1}{2}\right)\pi^{-(s+1)/2}L(\chi,s)&
\text{if $\chi(-1)=-1$\,.}
\end{cases}\end{equation}
Using $\theta(\chi,z)$ and (\ref{theta5}) one shows (Th.\ 7.19 of \cite{kn1})
\begin{equation}\label{lserieschi2}
\Lambda(\chi,s)= \begin{cases} w(\chi,m)\Lambda(\bar{\chi},1-s)
&\text{if $\chi(-1)=1$}\\
-iw(\chi,m)\Lambda(\bar{\chi},1-s)& 
\text{if $\chi(-1)=-1$\,.}
\end{cases}
\end{equation}
When $\chi=\chi_0$, the functional equation
for $\Lambda(\chi,s)$ is the one given above for the Riemann zeta
function.

We can loosely characterize the results just presented 
by saying that Dirichlet characters
are ``modular'' in that they arise from modular forms (albeit of
fractional weight).

Starting with Yutaka Taniyama in the 1950's,
 mathematicians began to suspect that the
connection between $L$-series of  abelian ``arithmetic objects'' defined
over $\Q$, such as Dirichlet characters, and modular forms might also
extend to ``non-abelian objects'' such as elliptic curves over
$\Q$. That such a
connection should exist at all is, at first glance, very surprising.
Indeed, the space $\mathcal H$ {\it already} has a deep connection
with elliptic curves as every elliptic curve over
$\C$ is isomorphic, as a complex analytic space, to
$E_z:=\C/L_z$ where $L_z$ is the lattice 
generated by $\{1,z\}$ for some
$z\in \mathcal H$. This new relationship between
$\mathcal H$ and elliptic curves, via $L$-series and modular forms,
is of a very different, and deeper, nature.

The modern rational for the existence of
this new connection is part of the general ``Langlands philosophy.'' 
To such an elliptic curve $E$ one associates an $L$-series (the definition
will be recalled below) $L(E,s)=L(E_\Q,s)$.
This $L$-series, and it twists by Dirichlet characters,
(also recalled below)
are conjectured to satisfy certain functional equations; in turn these
functional equations guarantee (due to Weil \cite{we1})
that the $L$-series arises from a modular
form of a specific type in essentially the same fashion as
$\zeta (s)$ arises from $\theta(z)$. 
(The existence of such functional equations and modular forms is
now, of course, well established, see below.)

As these ideas are crucial for us here, 
we will briefly recall them and refer the reader to \cite{kn1}
(for instance) for more details.
We begin by presenting more of the theory of modular forms. 
Let $\Gamma=SL_2(\Z)$, as above, and let $N$ be a positive integer. There is
clearly a homomorphism from $SL_2(\Z)$ to $SL_2(\Z/N)$ given by 
reducing the coefficients of the matrix modulo $N$. The kernel of this
mapping is denoted by $\Gamma(N)$. Any subgroup of $SL_2(\Z)$ which
contains $\Gamma (N)$, for some $N\geq 1$, is called a ``congruence
subgroup;'' we extend this notion to automorphisms of $\mathcal H$
in the obvious sense. (For instance, the group $\Gamma_\theta$
is a congruence subgroup in this sense, see \S I.5 of \cite{gu1}). 

Let $\tilde{\Gamma}$ be a congruence subgroup of $\Gamma$. From
now on we shall only consider modular forms for $\tilde{\Gamma}$
in the sense of
Definition \ref{unrestmod}; that is, without multipliers and where the
weight, $k$, is an integer.
The quotient space $\tilde{\Gamma}\backslash \mathcal H$ 
is an open Riemann surface that may be compactified
by adding a finite number of points called ``cusps;'' these cusps are in
one to one correspondence with $\tilde{\Gamma}\backslash
{\mathbb P}^1(\Q)$. For instance,
$\infty$ represents a cusp. The subgroup $\tilde{\Gamma}_\infty$ of
$\tilde{\Gamma}$ which fixes $\infty$ is of the form
$z\mapsto z+j$ where $j\in ({\mathfrak n})\subseteq \Z$ is an ideal
(and ${\mathfrak n}\geq 1$).
If $f(z)$ is a modular form for $\tilde \Gamma$, it then automatically
has a Fourier expansion 
\begin{equation}\label{expansion1}
f(z)=\sum_{n=-\infty}^\infty a_ne^{2\pi inz /\mathfrak n}
=\sum_{n=-\infty}^\infty a_n q_{\mathfrak n}^n\,,
\end{equation}
where $q_{\mathfrak n}=e^{2\pi i z/\mathfrak n}$.
Similar expansions are obtained at the other cusps by moving the cusp
to $\infty$ via an element of $SL_2(\Z)$. One says that the modular form
$f(z)$ is {\it holomorphic} if, at every cusp, all terms associated to
negative $n$ in the associated expansion (\ref{expansion1}) vanish. One
says that a holomorphic form is a {\it cusp form} if all terms associated to
$n=0$ at the cusps also vanish. 

One has holomorphic forms only if the weight $k$ is non-negative. Of course,
both the holomorphic and cusp forms of a given weight form vector spaces
of $\C$ which can be shown to be finite dimensional via standard
results on algebraic curves.

Of primary importance in the theory are the congruence subgroups
$\Gamma_0(N)\subset \Gamma$ defined by
\begin{equation}\label{gamma0def}
\Gamma_0(N):=\left\{ \left(\begin{array}{cc}a&b\\c&d\end{array}\right)
\in \Gamma\mid c\equiv 0\pmod{N}\right\}\,.\end{equation}
It is clear that for such subgroups one can choose
$\mathfrak n$ in the expansion (\ref{expansion1}) at $\infty$ equal to 
$1$, in which case we simply set $q=q_1=e^{2\pi i z}$. 

Now let $f(z)=\sum_{n=1}^\infty a_n q^n$ be a cusp form of weight $k$ for
$\Gamma_0(N)$. One sets 
\begin{equation}\label{Lcusp0}
L(f,s):=\sum_{n=1}^\infty a_n n^{-s}\,.
\end{equation}
In a fashion quite similar
to that of Equation (\ref{zeta1}), one finds
\begin{equation}\label{cusp1}
(2\pi)^{-s}\Gamma(s)L(f,s)=\int_0^\infty f(it)t^s\,\frac{dt}{t}\,.
\end{equation}
\noindent
Recall that the functional equation for $L(\chi,s)$ arises from the action 
$z\mapsto -1/z$ on $\mathcal H$. Similarly the functional equation for
$L(f,s)$ will arise from the action $z\mapsto -1/Nz$ on $\mathcal H$.
The matrix $\omega_N:=\left(\begin{array}{cc}0&-1\\N&0\end{array}\right)$
is not in $\Gamma_0(N)$ but rather in its normalizer. As such there is
a action of $\omega_N$ on cusp forms of a given weight for 
$\Gamma_0(N)$ and, to get a functional equation, one needs to further
assume that $f(z)$ is an eigenfunction for this action. The eigenvalue
$\varepsilon$ will be $\pm 1$. With this added assumption,
put
\begin{equation}\label{Lcusp2}
\Lambda(f,s):=N^{s/2}(2\pi)^{-s}\Gamma(s)L(f,s)\,.
\end{equation}

More generally let $\chi$ be a character (as above) which we now
assume has conductor $m$ which is prime to $N$. Set
\begin{equation}\label{twistedLcusp}
L(f,\chi,s):=\sum_{n=1}^\infty \chi(n)a_nn^{-s}\,,
\end{equation}
and
\begin{equation}\label{twistedLcusp2}
\Lambda(f,\chi,s):=(m^2N)^{s/2}(2\pi)^{-s}\Gamma(s)L(f,\chi,s)\,.
\end{equation}
One then has the functional equations
\begin{equation}\label{twistedLcusp3}
\Lambda(f,s)=\varepsilon (-1)^{k/2}\Lambda(f,k-s)\,,
\end{equation}
and
\begin{equation}\label{twistedLcusp4}
\Lambda(f,\chi,s)=\varepsilon (-1)^{k/2}w(\chi,m)\chi(-N)\Lambda(f,\bar{\chi},
k-s)\end{equation}
with $|w(\chi,m)|=1$.

We refer the reader to \cite{kn1} for a discussion of these functional 
equations which are due to Hecke 
(and which turn out, after all, to be very much in the
spirit of Riemann's theory for $\zeta(s)$ as in Equations (\ref{zeta2})
and (\ref{zeta3})).

Hecke also had a procedure for selecting those cusp forms
$f(z)$ for which $L(f,s)$ has an infinite product 
expansion (``Euler product'') similar to those
given in Equations (\ref{zeta0}) and (\ref{dirichlet1}) {\it except}
that the local factors will be degree $2$ polynomials in $p^{-s}$ for
almost all primes $p$. 
Hecke's idea can be very roughly sketched as follows: As above, every
point $z\in \mathcal H$ gives rise to
the elliptic curve $E_z:=\C/L_z=\C/\{\Z+\Z z\}$; therefore, one can
view modular forms  as certain 
functions on elliptic curves (together with, in the case of $\Gamma_0(N)$,
a cyclic subgroup $C$ of order $N$; for our purposes here we will simply
ignore such subgroups altogether). 
Given an integer $n$, we can associate to $E_z$ the set of all sub-lattices
$\{L_z^{(i)}\}$ of $L_z$ of
index $n$ as well as their associated elliptic curves
$\{E_z^{(i)}\}$. This association
depends only on $L_z$ and thus $E_z$. If $f$ is a function on
elliptic curves, we can then define a new function $T(n)f$ on elliptic
curves simply
by summing up the values of $f$ on the elliptic curves associated to
the sub-lattices; i.e.,
\begin{equation}\label{hecke1}
T(n)f(E_z):=\sum_i f(E_z^{(i)})\,.
\end{equation}
\noindent
It turns out that the ``Hecke operator,'' $f\mapsto T(n)f$, gives rise to linear
endomorphisms of both the space of modular forms and the space of cusp
forms of a given weight. These operators form a commutative
ring where $T(nm)=T(n)T(m)$ for coprime $n$ and $m$ but where,
for $p$ not dividing $ N$, $T(p^2)\neq T(p)^2$. The important
point is that those $L$-series, $L(f,s)$, which have
Euler products are precisely those which are associated to {\it eigenvectors}
(``eigenfunctions'') $f(z)$ for {\it all} Hecke operators $T(n)$. 
In order to establish this equivalence, one first shows that the 
Hecke eigenvalues are precisely the $q$-expansion coefficients of the
(normalized) form $f(z)$.

\subsection{$L$-functions of elliptic curves and elliptic modularity over $\Q$}\label{lseriesellcurve}
Now we can turn to elliptic curves over $\Q$. Such a curve $E$ is given
by a Weierstrass equation of the form 
\begin{equation}\label{weieqn}
y^2=x^3+ax+b
\end{equation}
 where $\{a,b\}\subset \Q$ and
$\Delta:=4a^3+27b^2$ is non-zero. The association
\begin{equation}\label{delta1}
z\in {\mathcal H}\mapsto \Delta\left(\C/L_z\right)
\end{equation}
makes $\Delta$ a cusp form of weight $12$ associated to the full modular
group $\Gamma$.

In order to discuss to define the local $L$-factors of $E$ we need to
discuss its reduction at the finite primes $p$ of $\Q$ following
\S VII of \cite{si1}; for more information we refer the reader there.
For simplicity
we begin by assuming that $p\geq 5$ with associated additive valuation
$v_p$. By the change of variables $(x,y)=(u^2\tilde{x},u^3\tilde{y})$,
$u\neq 0$, the Weierstrass equation 
(\ref{weieqn}) is changed into the Weierstrass
equation 
\begin{equation}\label{2weieqn}
\tilde{y}^2=\tilde{x}^3+\tilde{a}\tilde{x}+\tilde{b}\,,
\end{equation}
with $\tilde{a}=a/u^4$, $\tilde{b}=b/u^6$ and $\tilde{\Delta}=\Delta/u^{12}$.
Thus, by the appropriate choice of $u$, one can find a Weierstrass
equation for $E$ where all the coefficients are integral at $p$; in particular,
of course, $\Delta$ then is integral at $p$ also. Among all such equations,
the ones where $v_p(\Delta)$ is a minimum are called ``minimal
Weierstrass equations'' for $E$ at $p$. Such an equation is not unique
 but it is
easy to see that any two such equations will give rise to {\it isomorphic}
curves upon reduction modulo $p$ (obtained by reducing the coefficients of the
Weierstrass equation modulo $p$).

For almost all primes $p\geq 5$ (the ``good primes'')
the reduced minimal Weierstrass equation will also be an elliptic curve $E_p$
over $\Fp$. Let $n_p$ be the number of points on $E_p$ over $\Fp$ (which
is obviously an isomorphism invariant) and put
$a_p:=p+1-n_p$. Finally we define the local $L$-factor $L_p(E,u)$ by
\begin{equation}\label{ellipticfactor}
L_p(E,u):=\frac{1}{1-a_pu+pu^2}\,.
\end{equation}
A basic result, due to Hasse (Th.\ 10.5 of \cite{kn1}), establishes that 
$\dis L_p^{-1}(E,u)=(1-\alpha u)(1-\beta u)$ where $|\alpha|=|\beta|=p^{1/2}$.

Suppose now $E$ has ``bad'' reduction at a prime $p\geq 5$. Then from the 
Weierstrass equation one can see that the reduced
curve $E_p$ at $p$ must have either a node or a cusp. 
If the reduced curve has a node with slopes in $\Fp$ we say
that $E$ has ``split multiplicative reduction at $p$,'' if the reduced
curve has a node but where the slopes are not in $\Fp$ then we say
that $E$ has  ``non-split multiplicative reduction.'' If the reduced curve has
a cusp, then we say that {$E$ has ``additive reduction at $p$}.''
We can then describe the local factor at these finitely many bad
primes as follows.

\begin{equation}\label{ellipticfactor2}
L_p(E,u):=\begin{cases} \frac{1}{1-u}&\text{if $E$ has split multiplicative reduction at $p$}\\
\frac{1}{1+u}&\text{if $E$ has non-split multiplicative reduction at $p$}\\
1&\text{if $E$ has additive reduction at $p$\,.}\end{cases}\end{equation}

For the primes $p=2,3$ one has an exactly similar story {\it but} where
one has to use a more general form of the Weierstrass equation (\ref{weieqn});
again we refer the interested reader to \cite{si1} for the details.

\begin{rem}\label{additive}
We note for future use that if $E$ does not have good reduction at $p$, but
does acquire it over a finite extension, then $E_p$ must have
a cusp (since multiplicative reduction remains multiplicative reduction
over any finite extension; see Prop.\ 5.4.b of \cite{si1}). In this case we see from
Equation (\ref{ellipticfactor2}) that the
local factor $L_p(u)$ is identically $1$\end{rem}

Let $B$ be the finite set of bad primes for $E$. The {\it conductor
of $E/\Q$}, $N_E$, is defined by
\begin{equation}\label{conductor}
N_E:=\prod_{p\in B}p^{e_p}\end{equation}
where $e_p=1$ if $E$ has multiplicative reduction at $p$ and
$e_p\geq 2$ otherwise (in fact, equal to $2$ if $p\geq 5$);
see e.g., \S A.16 of \cite{si1}. The conductor
is a measure of how ``twisted'' the reduction of $E$ at bad primes 
actually is.

The $L$-series of the elliptic curve $E$, $L(E,s)$, is then defined
as 
\begin{equation}\label{lseriesellipticcurve}
L(E,s):=\prod_{\text{all primes $p$}} L_p(E,p^{-s})\,.
\end{equation}
Upon expanding the Euler product for $L(E,s)$ one obtains
\begin{equation}\label{lseriesellipticcurve2}
L(E,s)=\sum_{n=1}^\infty c_nn^{-s}\,.
\end{equation}
Let $\chi$ be a character of conductor $m$ prime
to $N:=N_E$. We then define the twisted $L$-series
$L(E,\chi,s)$ by
\begin{equation}\label{lseriesellipticcurve3}
L(E,\chi,s):=\sum_{n=1}^\infty\chi(n)c_nn^{-s}\,.
\end{equation}
One puts
\begin{equation}\label{lseriesellipticcurve4}
\Lambda(E,s):=N^{s/2}(2\pi)^{-s}\Gamma(s)L(E,s)\,,
\end{equation}
and 
\begin{equation}\label{lseriesellipticcurve5}
\Lambda(E,\chi,s):=(m^2N)^{s/2}(2\pi)^{-s}\Gamma(s)L(E,\chi,s)\,.
\end{equation}
It was conjectured (and is now a theorem) that $\Lambda(E,s)$
satisfies a functional equation of the form
\begin{equation}\label{lseriesellipticcurve6}
\Lambda(E,2-s)=\pm\Lambda(E,s)\,.\end{equation}
The sign $\pm 1$ here can expressed
as a product over all places of $\Q$ of ``local signs.'' Moreover
the sign of $E$ at $\infty$ is $-1$ and the sign at all good
primes is $+1$. In particular, it is remarkable that the sign is
then completely determined by the local signs at the {\it bad}
primes; see e.g., \cite{roh1},\cite{ha1}, \cite{ko1} and \cite{ri1}.
Similarly, when the conductor of $\chi$ is prime
to $N$, $\Lambda(E,\chi,s)$ was conjectured to satisfy (and is now
known to satisfy)
\begin{equation}\label{lseriesellipticcurve7}
\Lambda(E,\chi,s)=\pm w(m,\chi)\chi(-N)\Lambda(E,\bar{\chi},2-s)\,.
\end{equation}

As one can see, the functional equations (\ref{lseriesellipticcurve6})
and (\ref{lseriesellipticcurve7})
are remarkably like the functional equations 
given above (in (\ref{twistedLcusp3}) and (\ref{twistedLcusp4}))
for $L(f,s)$
and $L(f,\chi,s)$ where $f$ is a cusp form of weight $2$.
This ultimately led to the amazing expectation 
(the ``Modularity Conjecture'') that for
every $E$ one could find a cusp form $f_E(z)$ of weight $2$ for
$\Gamma_0(N)$
such that $f_E(z)$ is an eigenform for all the Hecke operators and
$L(E,\chi,s)=L(f_E,\chi,s)$ for all $\chi$ (of conductor
prime to $N$). In particular,
the conjectured analytic properties of $L(E,\chi,s)$ then follow
immediately from those of $L(f_E,\chi,s)$. Indeed, the results of
Weil \cite{we1} characterizes those Dirichlet series which
arise from cusp forms (for $\Gamma_0(N)$) as {\it precisely}
those satisfying functional equations (\ref{lseriesellipticcurve6})
and (\ref{lseriesellipticcurve7}).

The local $L$-factors of $E$ can also be obtained from Galois representations
associated to the elliptic curve as we will now explain.
Let $\ell$ be a prime number. Then to each
$\ell$ one attaches to $E$
the ``$\ell$-adic Tate module $T_\ell(E)$,''  defined as
the inverse limit of the $\ell^n$-division points on $E$ (\S III.7 of 
\cite{si1}).
One sees that the Tate module is a free $\Z_\ell$-module of rank $2$, 
and it's cohomological dual is defined by
\begin{equation}\label{h1}
H^1(E,\Q_\ell):=\text{Hom}_{\Z_\ell}(T_\ell(E),\Q_\ell)\,.
\end{equation}
Both $T_\ell(E)$ and $H^1(E,\Q_\ell)$ are naturally modules for
the Galois group $G:=\text{Gal}(\bar{\Q}/\Q)$ (where $\bar\Q$ is a chosen
algebraic closure), and one sees readily that this module is non-abelian
(in that the Galois action factors through a non-abelian group). 

The $L$-series $L(E,s)$ can also be expressed in terms of this ``compatible
system'' of representations on $H^1(E,\Q_\ell)$ 
(for varying $\ell$). At the good primes
one obtains $L_p(u)^{-1}$ as the characteristic polynomial of the geometric
Frobenius automorphism and at the bad primes one must first take
the fixed subspace of the inertia elements (and then use a geometric Frobenius
element etc.) see \cite{ta1} . 
This construction is the non-abelian version of the
procedure used to define $L(\chi,s)$ for Dirichlet characters. It also
justifies viewing the Modularity Conjecture as a non-abelian extension
of the relationship between characters and theta-functions given
in Subsection \ref{thetadir} as class field theory equates abelian
characters of $G$ with Dirichlet characters.

Once one knows that $L(E,s)=L(f_E,s)$ for a cusp form $f_E$ of weight $2$
(or, more technically correct, a ``newform'')
associated to $\Gamma_0(N)$ for some $N$, another remarkable sequence
of results take over. Indeed, as mentioned above,
the quotient space $\Gamma\backslash \mathcal H$ is compactified by
attaching cusps and can be realized as a smooth projective curve
$\bar{X}_0(N)$ over $\Q$. Using $f_E$ one  constructs (Th.\ 11.74 of
\cite{kn1})
a certain elliptic curve $E^\prime/\Q$ which is a quotient of the Jacobian of 
$\bar{X}_0(N)$. By construction one also 
has $L(E^\prime,s)=L(f_E,s)$. The existence and properties of 
of $E^\prime$ are due to Eichler and Shimura, as is the identification
of $L(E^\prime,s)$ with $L(f_E,s)$ (via an ``Eichler-Shimura relation'' which
connects the action of $T(p)$, $p$ prime, with the Frobenius automorphism at
$p$).

Clearly one should expect some
relationship between $E$ and $E^\prime$. 
Indeed, recall that an ``isogeny'' between
elliptic curves $E_1$ and $E_2$ is a surjective
map (of elliptic curves) $E_1\to E_2$
and once one has an isogeny $E_1\to E_2$ it
is easy to construct an isogeny $E_2 \to E_1$. If 
both $E_1$ and $E_2$ are defined over $\Q$,
and if the map between them is {\it also} defined over $\Q$, then
both elliptic curves will have the same local factors and $L$-series. 
Faltings \cite{fa1} tells us that the converse is also true; i.e., two
elliptic curves over $\Q$ with the same $L$-series (and thus
the same local factors) are then automatically
related by an isogeny (or are ``isogenous''). In particular, we conclude
that $E$ and $E^\prime$ are isogenous. Consequently, using modular
forms and curves, one obtains an amazing dictionary of the isogeny
classes of elliptic curves over $\Q$!

As we have stated, the modularity conjecture for elliptic curves
over $\Q$ is now a theorem.
The main work in establishing it was due to A.\ Wiles and 
Wiles and R.\ Taylor, \cite{wi1}, \cite{tawi1}.
The proof was then finished in \cite{di1}, \cite{cdt1}, and \cite{bcdt1}.

A key part of the Wiles' proof is the result of R.\ Langlands
and J.\ Tunnell, \cite{tu1}, \cite{lan1}.  This establishes that certain
non-abelian and  
complex-valued (and thus of finite image)
representations of the absolute Galois group
of $\Q$ are modular in a similar 
sense as Dirichlet characters (in fact, one may view these
representations as being non-abelian generalizations of Dirichlet
characters).
In other words the $L$-functions associated to these representations
arise from certain cusp forms, which, in this case, are of weight $1$.
This gives yet another instance of the crucial
role played by modular forms in classical arithmetic.

\section{Elliptic modularity for $k=\Fr(T)$}\label{ellmodffields}
In this section we will explain how the classification of isogeny
classes of elliptic curves over $\Q$ can be translated to the case
of the global function field $k=\Fr(T)$ for a certain
class of elliptic curves over $k$ (e.g., those elliptic curves with
split-multiplicative reduction at the place $\infty$ defined below). 
For a clear and thorough
exposition of these ideas, we refer the reader to \cite{gr1}.

\subsection{The $L$-series of an elliptic curve over $k$}\label{Lseriescurvek}
Let $\Fr$ be the finite field with $r=p^m$ elements and $p$ prime. 
We let $k=\Fr(T)$ be the rational field in an indeterminate $T$. 
Let $E$ be an elliptic curve over $k$.  The $L$-series of $E/k$ is defined
in a completely analogous fashion to that of $E/\Q$ as given above.
More precisely, let $w$ be a place of $k$ with local ring
${\mathcal O}_w$ and associated finite field $\F_w$.
Put $Nw=|\F_w|$, which is a power of $r$. If $w$ is a place of
good reduction then one sets $n_w$ to be the number of points on the reduction
$E_w$ and $a_w:= Nw+1-n_w$.The local $L$-factor is
\begin{equation}\label{localLfactorffield}
L_w(u):=\frac{1}{1-a_wu+(Nw) u^2}\,.\end{equation}
If $w$ is a place of bad reduction, then the local factor is defined
as in (\ref{ellipticfactor2}). Finally, we put
\begin{equation}\label{Lseriesellipffield}
L(E,s)=L(E/k,s):=\prod_{\text{all $w$}}L_w(Nw^{-s})\,.\end{equation}
\subsection{Automorphic representations}
In Subsections \ref{thetadir} and
\ref{modularforms}, we sketched some important aspects of the theory
of classical elliptic modular forms. Using these forms, one can obtain
the analytic properties of various $L$-series of arithmetic
objects over $\Q$ as described above. However, the general formalism and
construction of $L$-series can be given in far greater generality than 
just for objects over $\Q$. For instance, one can work with arbitrary
global fields (such as $k$ in the last subsection). 

The correct generalization of cusp forms that works for
arbitrary global fields is the notion of a
``cuspidal automorphic representation'' (see, e.g., \cite{jl1}, 
\cite{bj1}, \cite{bor1}, \cite{lan2}). For our
purposes, it is sufficient to view such automorphic representations
as being ``generalized cuspidal (Hecke) eigenforms.'' 
Indeed, a cuspidal automorphic
representation $\pi$ can be given an associated $L$-series $L(\pi,s)$ 
which arises from an Euler product and which has
good properties (such as functional equations). Moreover, the
$L$-series of an arithmetic object is {\it always} conjectured
in the Langlands philosophy to equal the
$L$-series of a certain associated cuspidal automorphic representation $\pi_E$
(which completely generalizes the association $E\leftrightarrow f_E$
discussed in Subsection \ref{lseriesellcurve}).
In the function field case, this is now known to
be a theorem to due to the labors of Drinfeld, Lafforgue \cite{laf1} 
and others,
(an excellent source in this regard is \cite{lau1} and its references).

In particular, due to the cohomological results of Grothendieck \cite{gro1}
one knows that the $L$-series
$L(E,s)$, and its twists by abelian characters (the generalization of
$L(E,\chi,s)$, see (\ref{lseriesellipticcurve5})) have functional equations
(in fact, one knows that, in this case, they are polynomials in
$r^{-s}$).
Thus the general theory of automorphic representations will associate to
$E/k$ an automorphic representation $\pi_E=\pi_{E/k}$ with 
$L(\pi_E,s)=L(E,s)$ \cite{de1}. The reader
should realize that this is very different than the case of elliptic
curves $E$ over $\Q$ where one deduces the analytic properties of
$L(E,s)$ at the same time as one finds the associated modular form
$f$. In the case of $E/k$ we know the analytic properties of the 
$L$-series directly,
without having the associated $\pi_{E/k}$; in fact, one constructs
$\pi_{E/k}$ from this knowledge.

What is lacking in the function field case is a concrete realization
of the isogeny class of $E$, as was accomplished in the case of elliptic
curves over $\Q$ via the Jacobians of elliptic modular curves. It is precisely
here that the work of V.G.\ Drinfeld comes in.

\subsection{A Quick introduction to Drinfeld modules}\label{drinfeld}
In 1973, V.G.\ Drinfeld introduced his ``elliptic modules'' \cite{dr1} which
are now called ``Drinfeld modules'' in his honor. The 
analytic construction of Drinfeld modules is based on that of elliptic curves
where the Archimedean place is singled out. Thus one begins by singling 
out a particular ``infinite'' 
place of $k=\Fr(T)$; the obvious one is ``$\infty$'' where
$v_\infty(1/T)=1$ (so, naturally, $\infty$ corresponds to the usual point
$\infty\in {\mathbb P}^1(\Fr)$). The ring $A:=\Fr[T]$ consists of those
rational functions which are regular away from $\infty$. The field
$K:=k_\infty=\Fr((1/T))$ is a local field which contains $A$ discretely
and $K/A$ is compact. The standard analogy is with $\Z\subset \R$
with $\R/\Z$ being compact. The algebraic closure of $K$, denoted
$\bar{K}$, is infinite dimensional over $K$ and is not
complete. However, $v_\infty$ 
lifts to $\bar{K}$ in a canonical way and every subextension
$L\subset \bar{K}$ which is finite dimensional over $K$ is, in fact,
complete. Thus we may use analytic methods over $L$. 

A $\Z$-lattice in $\C$ is a discrete (in the standard topology on
$\C$) $\Z$-submodule which may then be shown to 
have rank $1$ or rank $2$ (reflecting the fact that $[\C \colon \R]=2$).
The rank two lattices are precisely those that give rise to elliptic curves.
An {\it $A$-lattice} $M$ is a finitely generated, discrete (i.e.,
finitely many elements in any bounded ball with the metric
generated by $v_\infty$), $A$-submodule of $\bar{K}$. As $M$ is finitely
generated and obviously torsion-free, it is free of some rank $t=t_M$
and generates a finite extension of $K$.
As $[\bar{K}\colon K]=
\infty$, one can have lattices of {\it arbitrary} rank. 
To $M$ one attaches its exponential function
\begin{equation}\label{expfunction}
e_M(x)=x\prod_{0\neq m\in M}(1-x/m)\,.\end{equation}
As $M$ is discrete, it is easy to see that $e_M(x)$ converges for
all $x\in \bar K$; that is, $e_M(x)$ is an entire non-Archimedean
function.
As $M$ is finitely generated, the Taylor coefficients of $e_M(x)$ will
lie in some finite extension of $K$. Consequently, if $x\in \bar K$ then
$e_M(x)$ converges to an element of $\bar K$.

Non-Archimedean analysis is highly algebraic in nature. In particular,
like polynomials, {\it all} entire non-Archimedean functions in
$1$-variable are surjective
(as a function on $\bar K$) and are 
determined up to a constant by their divisors.

The main ``miracle'' of $e_M(x)$ is that the map $e_M(x)\colon \bar{K}\to
\bar{K}$ is actually $\Fr$-linear; thus one has $e_M(x+y)=e_M(x)+e_M(y)$.
This is due to Drinfeld \cite{dr1} (but uses
some combinatorial arguments on polynomials that have been known for ages).
One then deduces the remarkable isomorphism of $\Fr$-vector spaces
\begin{equation}\label{expiso}
e_M(x)\colon \bar{K}/M \to \bar{K}\,.
\end{equation}
The idea behind the analytic construction of Drinfeld modules is to carry
over the natural quotient 
$A$-module structure on the left of (\ref{expiso}) to 
$\bar{K}$ via $e_M(x)$ (just as one carries over the $\Z$-module structure
on $\C/\{\Z+\Z z\}$ to the associated elliptic curve). 

More precisely, let $a\in A$ be a polynomial of degree $d$ and let
\begin{equation}\label{Ea}
E_a:=\{e_M(\alpha)\mid  \alpha\in a^{-1}M/M\}\,.\end{equation}
Thus $E_a\subset \bar{K}$ is a vector space over $\Fr$ of dimension $dt_M$ (and so
$|E_a|=r^{dt_M}$). Put,
\begin{equation}\label{drinlat}
\psi_a(x)=ax\prod_{0\neq \alpha\in E_a}(1-x/\alpha)\,.\end{equation}
As $E_a$ is a finite set, $\psi_a(x)$ is a {\it polynomial}. The same
combinatorics as used for $e_M(x)$ also establishes that $\psi_a(x)$ is an
$\Fr$-linear function. As non-Archimedean entire functions are determined
up to a constant by their divisors, a simple calculation then gives the
basic identity
\begin{equation}\label{basicdrinfeld}
e_M(ax)=\psi_a(e_M(x))\,.\end{equation}
That is, the standard $A$-action on $\bar{K}$,
$\{a,x\}\mapsto ax$ (on the left hand side of
(\ref{expiso})) gets transfered over to the action
$\{a,x\}\mapsto \psi_a(x)$ (on the
right hand side of (\ref{expiso})). 
In particular, $\bar K$ inherits a new $A$-module action
called a ``Drinfeld module.'' 

The mapping $a\in A\mapsto \psi_a(x)$ is readily seen to give an injection
of $A$ into the algebra of $\Fr$-linear polynomials (with composition of
polynomials as multiplication). 
As it is an algebra map, it is uniquely determined by
$\psi_T$ which, as an $\Fr$-linear polynomial, is given by
\begin{equation} \label{psiT}
\psi_T(x)=Tx+\sum_{i=1}^{t_M}a_ix^{r^i}\,,\end{equation}
where $a_{t_M}\neq 0$.
The rank of the lattice $M$, $t_M$, is also the ``rank'' of the
Drinfeld module. 

As $\psi_a(x)$ is a polynomial in $x$ for all $a\in A$, the notion of
a Drinfeld module is really an algebraic one exactly as is the case
with elliptic curves. Thus it makes sense over any field
$L$ containing $A/\mathfrak p$ for any prime ideal $\mathfrak p$
of $A$ (including, obviously, ${\mathfrak p}=(0)$). Indeed, 
to get a rank $t$ Drinfeld $A$-module over $ L$, for a positive integer
$t$, as done in (\ref{psiT}) one just needs $\{a_i\}_{i=1}^t\subseteq L$ with
$a_t\neq 0$ and ``$T$'' represents its image in $L$. It is common to
denote this image by ``$\theta$''; the use of $\theta$ allows us to distinguish
when ``$T$'' is an operator via a Drinfeld module and ``$T=\theta \in L$'' is
a scalar.

Any Drinfeld module of rank $t$ defined over $\bar K$ can be shown to
arise from a lattice $M$ of the same rank via $e_M(x)$ as above;
this is in exact agreement with the analytic theory of elliptic curves.
For more on Drinfeld modules we refer the reader to  \cite{hay1}
or \cite{go4}.

\begin{example}\label{carlitz1}
Let $C$ be the rank $1$ Drinfeld module over $\Fr(T)=\Fr(\theta)$ defined by
\begin{equation}\label{carlitz2}
C_T(x):=Tx+x^r=\theta x+x^r\,.\end{equation}
It is clear that $C$ has rank $1$ and is the simplest possible
Drinfeld module. It is called the ``Carlitz module'' after the
work of L.\ Carlitz \cite{ca1}. It is associated to a rank one
lattice $M:=A\xi$ where $0\neq \xi\in \bar{K}$ and $\xi^{r-1}\in K$.
\end{example}

\begin{rem}\label{manycases}
As mentioned in the
introduction, the theory of Drinfeld modules actually exists in much greater 
generality where one replaces $\Fr(T)$ by an arbitrary
global function field $k$ 
of characteristic $p$ and $\Fr[T]$ by the affine algebra $A$
of functions regular away from a fixed place ``$\infty$'' of $k$.
The ring $A$ is readily seen to be a Dedekind domain with finite
class group and unit group equal to $\Fr^\ast$.
\end{rem}

\subsection{The rigid space $\Omega^2$}\label{omega}
A Drinfeld module $\psi$ of rank $2$ over $\bar{K}$ is given by
\begin{equation}\label{rank2drinfeld}
\psi_T(x)=Tx+g(\psi)x^r+\Delta(\psi) x^{r^2}=\theta x+g(\psi)x^r
+\Delta(\psi)x^{r^2}\,,\end{equation}
where $\{g(\psi),\Delta(\psi)\}\subset \bar K$ and 
$\Delta(\psi)\neq 0$. From our last subsection, we 
know that $\psi$ arises from a rank $2$ 
$A$-lattice $M$ of the 
form $Az_1+Az_2$ where the discreteness of $M$ is equivalent to
$z_1/z_2\in {\bar K}\backslash K$.

\begin{defn}\label{omega2}
We set $\Omega^2:={\bar K}\backslash K\,.$\end{defn}
The space $\Omega^2$ was defined by Drinfeld in \cite{dr1}
and, in fact, there is
an $\Omega^i$ for all $i=1,2\cdots$. As we are only interested
in $\Omega^2$ here, from now on we shall simply denote it ``$\Omega$.''
The space $\Omega$ is clearly analogous to the $\C\backslash \R$, which,
in turn, is precisely the upper and lower half-planes. Like
$\C\backslash \R$, $\Omega$ has an analytic structure which allows one to
use analytic continuation. This structure is called a ``rigid analytic
space.'' Surprisingly, with this rigid structure $\Omega$ becomes a connected
(but not simply connected) space unlike, of course,
the classical upper half-plane $\mathcal H$. Rigid analysis allows
one to handle non-Archimedean functions in a manner very similar to that
of complex analytic functions.

The space $\Omega$ has an action of $\Gamma:=GL_2 (A)$ on it
completely analogous to the classical action of $SL_2(\Z)$
on $\mathcal H$. Let
$\gamma=\left(\begin{array}{cc}a&b\\c&d\end{array}\right)\in GL_2(A)$ and
$z\in \Omega$. Then, exactly as in (\ref{gammaact}), we
set
$$\gamma z:=\frac{az+b}{cz+d}\,;$$
the map $z\mapsto \gamma z$ is a rigid analytic automorphism of
$\Omega$ where the inverse transformation arises from the inverse matrix.

Let $0\neq N\in A$. The definition of $\Gamma_0(N)\subseteq \Gamma$ is
exactly the same as in the classical case (\ref{gamma0def}). The rigid
analytic space $X_0(N):=\Gamma_0(N)\backslash \Omega$ may 
then be realized as the
underlying analytic space associated to an affine algebraic curve. 
As in the classical case, this space is compactified by adjoining
a finite number of ``cusps'' and these cusps are given by
$\Gamma_0(N)\backslash {\mathbb P}^1(k)$; we denote the compactified space
by $\bar{X}_0(N)$. As in the number field case, $\bar{X}_0(N)$ may
be realized canonically as a complete, smooth, geometrically
connected curve over $k$.
Analogs of the constructions of Eichler and Shimura for elliptic
curves can then be given on
the Jacobian of $\bar{X}_0(N)$.
\subsection{Modularity}\label{drinmod}
Finally, Drinfeld \cite{dr1} established a reciprocity
law which, in particular, identifies those 
automorphic representations given by elliptic curves over $k$ which
occur as quotients of the Jacobians of $\bar{X}_0(N)$.

Modularity for the class of elliptic curves over
$k$ with split-multiplicative
reduction at $\infty$ is then accomplished via Drinfeld's reciprocity
law coupled with the results of 
Grothendieck and Deligne mentioned above (which allow us to
construct an associated automorphic representation),
and the work of Y.\ Zarhin \cite{za1}, \cite{za2} ($p\neq 2$, the case
$p=2$ is unpublished work of S.\ Mori; see
Cor.\ XII.2.4 and Th.\ XII 2.5 of \cite{mb1}) establishing that
the local factors of the elliptic curve over $k$ determine
its isogeny class (as Faltings was later to show for number fields).

\section{Modularity for Drinfeld modules?}\label{moddrinfeld}
As we pointed out Subsection \ref{modularforms}, classical elliptic
modular curves play two (at least) distinct roles in number theory. The
first role is their use in classifying elliptic curves
(with various level structures) up to isomorphism. The second,
very different, role lies in their use classifying elliptic curves over 
$\Q$ up to isogeny.

Let $k=\Fr(T)$ as in the last section. Then we have seen how,
for elliptic curves over $k$ with split multiplicative
reduction at $\infty$, one classifies their isogeny classes 
via the moduli space of Drinfeld modules of rank $2$. One is thus
lead to ask whether Drinfeld modules themselves are
``modular'' in some reasonable sense. It is our goal here
to explain finally how this may indeed be possible. In particular, just as the
$L$-series of elliptic curves plays a crucial role in the modularity
conjecture for elliptic curves over $\Q$, so too will the $L$-series
of Drinfeld modules play an essential role here. 

The basic idea is, roughly, that a Drinfeld module (or related object) will
be called ``modular'' if its $L$-series can be obtained from the
$L$-series of a rigid analytic cusp form under a simple translation  of the
argument.

The definition of such $L$-series proceeds
very much like that the definition of
$L$-functions of elliptic curves; one defines first the local
Euler factor associated to a prime
and then takes their product. We elaborate this
construction first for the Carlitz module.

\subsection{The $L$-series of the Carlitz module over $k$}\label{Lcarlitz}
Let $w=(f)$ be the prime ideal associated to a monic prime $f\in A$. Let
$\F_w:=A/w$ be the associated finite field and let $C(\F_w)$
be $\F_w$ viewed as an $A$-module via the Carlitz action. As $C$ has
rank $1$, it is easy to see that $C(\F_w)$ is isomorphic to
$A/(g)$ for some monic $g\in A$. A simple calculation implies that
$C_f(x)\equiv x^{r^{\deg f}}\pmod{w}$; thus $f-1$ annihilates
$C(\F_w)$. Therefore $g\mid (f-1)$ and counting points implies
that, in fact, $g=f-1$. Consequently,
to $w$ and $C$ we define the local $L$-factor $L_w(C,u)$ by
\begin{equation}\label{charpL1}
L_w(C,u)=L_f(C,u):=\frac{1}{(1-fu)}\,,\end{equation}
which the reader will see is a rank $1$, characteristic $p$,
version of (\ref{ellipticfactor}).

We set
\begin{equation}\label{sinfty}
S_\infty:=\bar{K}^\ast\times\Zp\,,
\end{equation}
which is a topological
abelian group whose group operation will be written additively.
For $s=(x,y)\in S_\infty$ and $a\in A$ monic we define
\begin{equation}\label{as}
a^s:=x^{\deg a}\cdot (a/T^{\deg a})^y\,,\end{equation}
where $(a/T^{\deg a})^y$ is defined using the binomial theorem
(and converges in $K$ as $a/T^{\deg a}=1+$\{higher terms in $1/T$\}).
In particular, note that if $s_i:=(T^i,i)$, $i\in \Z$,
then $a^s=a^i$; as such, we shall
commonly write ``$i$'' for $s_i=(T^i,i)\in S_\infty$.
One views $S_\infty$ as a topological 
abelian group with
the integers embedded (discretely) as a subgroup.

We now define the $L$-function $L(C,s)$, $s\in S_\infty$,
of the Carlitz module by
\begin{equation}\label{Lcarlitz1}
L(C,s)=L(C,x,y):=\prod_\text{$f$ monic prime}L_f(C,f^{-s})=\prod_f (1-ff^{-s})^{-1}=\prod_f(1-f^{1-s})^{-1}\,.
\end{equation}
Upon expanding (\ref{Lcarlitz1}), we find
\begin{equation}\label{Lcarlitz2}
L(C,s)=\sum_\text{$n$ monic}n^{1-s}=\sum_{e=0}^\infty x^{-e}
\left(\sum_{\substack{\text{$n$ monic}\\ \deg n=e}}n\langle n\rangle^{-y}
\right)\,.\end{equation} 
In this case, elementary estimates (\S 8.8 of \cite{go4})
allow us to establish 
that $L(C,x,y)$ is an entire power series for all $y\in \Zp$.
Moreover, the resulting function on $S_\infty$ is also continuous 
and its ``zeroes flow continuously.'' (The best technical
definition of this concept is via non-Archimedean Frechet spaces
as in \cite{boc1}.)

Now let $y=-i$ for $i$ a non-negative integer. The same elementary
estimates also allow us to show that $L(C,x,-i)$ is a polynomial
in $x^{-1}$; one then immediately deduces that $L(C,x/T^i,-i)\in A[x^{-1}]$.
As the set of non-positive integers is dense in
$\Zp$, we see that the set of {\it special
polynomials} $\{L(C,x/T^i,-i)\}$ determines $L(C,s)$ as a function
on $S_\infty$. In Subsection \ref{crystals} we will see that such polynomials
are {\it cohomological} in nature which will be the key towards
handling the $L$-series of an arbitrary Drinfeld module.

\begin{rem}\label{carlitzspecialvalue}
Implicit in the definition of $L(C,s)$ is the ``zeta-function of $A$''
defined by 
\begin{equation}\label{zetafunction}
 \zeta_A(s)=\prod_{\text{$f$ monic prime}}(1-f^{-s})^{-1}\,;\end{equation}
so $L(C,s)=\zeta_A(s-1)=\zeta_A(s-s_1)$ ($s_1$ as above). Clearly the analytic
properties of $\zeta_A(s)$ follow from those of $L(C,s)$. Let
now $i=s_i$ be a positive integer which is divisible by $(r-1)$
and let $\xi$ be the period of the Carlitz module (as in Example \ref{carlitz2}).
It is then easy to see that $0\neq \zeta_A(i)/\xi^i\in \Fr(T)$ which is a version of
the classical result of Euler on zeta-values at positive
even integers.\end{rem}

\begin{rem}\label{vadic}
It is natural to wonder if there is some obstruction to interpolating
the set 
$\{L(C,x/T^i,-i)\}$ at a finite prime $v$ of $k$. In fact, there is none
(see \S 8 of \cite{go4} or \cite{boc1}).
Just as one obtains functions on $\bar{K}^\ast \times \Zp$, so too
does one obtain functions on $\bar{k}_v^\ast\times S_v$ where $\bar{k}_v$ is
the algebraic closure of the completion $k_v$ and 
$ S_v:=\varprojlim_j {\Z}/(p^j(r^{\deg v}-1))$. These functions
have remarkably similar properties to $L(C,s)$, $s\in S_\infty$. While
we do not emphasize such $v$-adic functions here for space considerations,
their existence is an important and natural part of the theory.
\end{rem}

\subsection{A quick introduction to
$T$-modules and $\tau$-sheaves}\label{taumodule}
In order to understand the general $L$-series of Drinfeld modules,
and their possible ``modular'' interpretation we need to expand the
category of objects under study. 

We will begin first with $T$-modules. This is an idea due to Greg
Anderson \cite{an1} (see also \S 5 of \cite{go4}), 
based on Drinfeld's notion of ``shtuka'' or
``elliptic sheaf.'' The idea behind it
is to replace the use of polynomials in $1$ variable in the definition of
a Drinfeld module in Subsection \ref{drinfeld} with polynomials in
many variables. 

Thus let $L$ be any extension of $\Fr$ and consider the algebraic
group $E:={\mathbb G}_a^e$ over $L$, where ${\mathbb G}_a$ is the additive
group. Let $\dis x=\left(\begin{array}{c} x_1\\\vdots\\x_e\end{array}\right)
\in E$. There are two obvious types of $\Fr$-linear endomorphisms
of $E$ as an algebraic group. The first is the $r^i$-th power mapping 
defined by $\dis x^{r^i}:=\left(\begin{array}{c} x_1^{r^i}\\\vdots\\
x_e^{r^i}\end{array}\right)$. The second is $x\mapsto Mx$ where $M
\in M_e(L)=$\{$e\times e$ matrices with coefficients in $L$\}. 
It is then easy to 
see that any $\Fr$-linear endomorphism of $E$ is just a combination of
these; i.e., it can be written $\sum_{j=0}^t M_ix^{r^i}$, for $M_i\in
M_e(L)$. We let $\text{End}_\Fr(E)$ be the set of all $\Fr$-linear
endomorphisms viewed as an $\Fr$-algebra under composition. 

As is standard, we let $I_e\in M_e(L)$ be the identity matrix.

Now let $L$ be a field equipped with an $\Fr$-algebra map $\imath\colon
A=\Fr[T]\to L$. We set $\theta:=\imath (T)$ as before.
A {\it $T$-module over $L$} is then a pair
$E=(E_\text{gp},\psi_E)$ where $E_\text{gp}$ is an algebraic group isomorphic to ${\mathbb G}_a^e$,
for some $e$,  and $\psi=\psi_E\colon A\to \text{End}_\Fr(E_\text{gp})$ is an injection
of $\Fr$-algebras. This injection is uniquely determined by $\psi_T$
which is further required to satisfy
\begin{equation}\label{Tmoduledef}
\psi_T(x)=(\theta I_e +N)x+\sum_{i=1}^t M_i x^{r^i}\,,
\end{equation}
for some (possibly vanishing) 
$M_i\in M_e(L)$ and {\it nilpotent} $N\in M_e(L)$. 

\begin{rem}\label{francis}
One can restate (\ref{Tmoduledef}) as 
\begin{equation}\label{francis2}
\psi_T(x)=\Theta x+\sum_{i=1}^t M_i x^{r^i}\,,\end{equation}
where $\theta$ is the only eigenvalue for $\Theta\in M_e(L)$ 
(i.e., the characteristic polynomial of $\Theta$ is
$(\lambda-\theta)^e$).\end{rem}

The reader may well wonder why one allows the existence of the
nilpotent matrix $N$ in (\ref{Tmoduledef}). The reason is that
it's existence allows us to introduce a {\it tensor product} into
the theory, \cite{an1}.

\begin{example}\label{trivialTmodule}
Let $e$ be arbitrary and set
$\psi_T(x)=\theta I_e x=\theta x$.
This is indeed a $T$-module under the above definition albeit a 
not very interesting one. Furthermore,
note that when $e=1$, we do {\it not} get
a Drinfeld module. Indeed, this ``trivial $T$-module'' is precisely the
case ruled out in the definition of Drinfeld modules.\end{example}

\begin{example}\label{carlitztensorpowers}
(See \cite{at1}.)
Let $L=k=\Fr(T)$ and let $\imath$ be the identity mapping. 
Let $n$ be a positive integer and $C^{\otimes n}_\text{gp}:={\mathbb G}_a^n$.
Let $N_n$ be the $n\times n$ matrix
$$\begin{pmatrix} 0 &1 &&0\\
&\ddots &\ddots&\cr
&&\ddots &1\cr
0 &&&0\end{pmatrix}\,,$$
and $V_n$ the $n\times n$ matrix
$$\begin{pmatrix}0 &\dots &0\cr
\vdots &&\vdots\cr
1 &\dots &0\end{pmatrix}\,.$$
We then set $C^{\otimes n}$ to be the injection of
$A$ into $\text{End}_\Fr(C^{\otimes n}_\text{gp})$ given by
\begin{equation}\label{carlitztensorpowers2}
C^{\otimes n}_T(x):=(\theta I_n+N_n)x+V_nx^r\,.\end{equation}
We then have a $T$-module ${\mathfrak C}^{\otimes n}:=(C^{\otimes n}_\text{gp},
C^{\otimes n})$.
We call ${\mathfrak C}^{\otimes n}$ the ``$n$-th tensor power
of the Carlitz module.'' 
\end{example}

One commonly uses $C^{\otimes n}$ interchangeably with 
${\mathfrak C}^{\otimes n}$. Clearly $C^{\otimes 1}$ coincides
with the Carlitz module $C$ as defined in Example \ref{carlitz1}.

In order to explain how the tensor product
appears in the theory, we begin with a dual construction
originally due to Drinfeld. Let $\psi$ be a Drinfeld module of rank
$d$ over a field $L$. Let
\begin{equation}\label{motive1}
M:=\text{Hom}^{(r)}_L({\mathbb G}_a, {\mathbb G}_a)
\end{equation}
be the vector space 
of $\Fr$-linear homomorphisms of the additive group to itself
as an algebraic group over $L$. We make the group $M$ into a left module
over $L\otimes_\Fr \Fr[T]\simeq L[T]$ via
$\psi$ as follows: Let $f(x)\in
M$, $a\in \Fr[T]$ and $l\in L$. Then we put
\begin{equation}\label{motive2}
l\otimes a\cdot f(x):=lf(\psi_a(x))\,.\end{equation}
It is easy to see (using a right division algorithm) that $M$ is
free over $L[T]$ of rank $d$. The module $M$ is called the ``motive''
of $\psi$. 

More generally, let $E=(E_\text{gp},\psi)$ be an arbitrary $T$-module. We define
its motive $M=M(E)$ as the group of $\Fr$-linear
morphisms of $E_\text{gp}$ to ${\mathbb G}_a$ over $L$, exactly as (\ref{motive1}).
The action of $L[T]$ on $M$ is defined as in (\ref{motive2}). 
The $T$-module $E$ is said to be {\it abelian} if and only if
its motive $M$ is finitely generated over $L[T]$. In this case,
$M$ is then free over $L[T]$ of finite rank which is also
the {\it rank of $E$}.

For instance, Drinfeld modules are exactly the
$1$-dimensional abelian $T$-modules. 
As an exercise, the reader may check
that $C^{\otimes n}$ of Example \ref{carlitztensorpowers}
is abelian (of rank $1$)
whereas the trivial $T$-module of Example \ref{trivialTmodule} is not.

The motive $M$ of a $T$-module also comes equipped with a canonical
endomorphism $\tau$ defined by
\begin{equation}\label{tauaction} \tau f(x):=f^r(x)\,.\end{equation}
Notice that $\tau (l f)=l^r\tau (f)$ for $l\in L\subset L[T]$ whereas
$\tau (a\cdot f)= a \cdot \tau(f)$ for $a\in A\subset L[T]$; we call such
a mapping ``partially Frobenius-linear.'' In Anderson's theory \cite{an1},
it is the interplay  between the $T$-action and the 
partially Frobenius-linear $\tau$-action that allows one
to pass back and forth between a $T$-module and its motive.

A ``$\tau$-sheaf'' is then just a globalization of $M$ viewed
as an $L[T]$-module equipped with the action of $\tau$. More precisely,
let $X$ be a scheme over $\Fr$. 
\begin{defn}\label{tausheaf} (See \cite{bp1} or \cite{boc1})
A {\it coherent $\tau$-sheaf} on $X$ is a pair $\underline{\mathcal F}:=
({\mathcal F},\tau)$ consisting of a coherent sheaf $\mathcal F$ on
$X\times_\Fr {\mathbb A}^1$ and a partially Frobenius-linear
mapping $\tau=\tau_{_{\mathcal F}}\colon {\mathcal F}\to \mathcal F$. A
{\it morphism} of $\tau$-sheaves is a morphism of the underlying
coherent sheaves which commutes with the $\tau$-actions.
\end{defn}
\noindent
Therefore
$M$, with the standard action of $\tau$ (\ref{tauaction}),
canonically  gives a $\tau$-sheaf on $\text{Spec}(L)\times {\mathbb A}^1$.
We call a $\tau$-sheaf $\underline{\mathcal F}$
{\it locally-free} if $\mathcal F$ is locally-free
on $X\times {\mathbb A}^1$. We call $\underline{\mathcal F}$ a {\it strict
$\tau$-sheaf} if it is locally-free and $\tau$ is injective. (Our strict
$\tau$-sheaves are the ``$\tau$-sheaves'' of \cite{ga1}.) 
The $\tau$-sheaves arising
from $T$-modules, for instance, are strict in this definition.

The {\it rank} of a strict $\tau$-sheaf is just the rank of the underlying
vector bundle.

\begin{example}\label{carlitztausheaf}
We will describe here the $\tau$-sheaf
$\underline{\mathcal C}=({\mathcal C},\tau)$ on
$\text{Spec}(\Fr(\theta))$
associated to the Carlitz module $C$.
The underlying space for the vector bundle is 
$\text{Spec}(\Fr(\theta))\times {\mathbb A}^1\simeq
\text{Spec}(\Fr(\theta)[T])$; for the moment let us call this
product $Y$. Over $Y$ the coherent module $\mathcal C$ given by 
$M=M(C)$ is isomorphic to the structure sheaf ${\mathcal O}_Y$.
The action of $\tau$ is then easily checked to be given by 
\begin{equation}\label{carlitztausheaf2}
\tau (\sum h_i(\theta)T^i):=(T-\theta)\sum h_i^r(\theta)T^i\,.\end{equation}
\end{example}

\begin{rem}\label{taumodulerankone}
Example \ref{carlitztausheaf} suggests the following general construction
of rank $1$ strict $\tau$-sheaves. Let $Y=\text{Spec}(\Fr(\theta)[T])$ as in
the example and let $g(\theta,T)$ be an arbitrary non-trivial function in
$\Fr(\theta)[T]$. We then define the $\tau$-sheaf $\underline{{\mathcal F}_g}$
to have underlying sheaf ${\mathcal O}_Y$ and $\tau=\tau_g$-action
given by
\begin{equation}\label{carlitztausheaf3}
\tau_g (\sum h_i(\theta) T^i):=g(\theta,T)\sum h_i^r(\theta)T^i\,.
\end{equation}
As as example, let $0\neq\beta \in \Fr(\theta)$. 
One then has the general rank $1$
Drinfeld module $C^{(\beta)}$ defined over $\Fr(\theta)$ by
\begin{equation}\label{taumodulerankone2}
C^{(\beta)}_T(x):=\theta x+\beta x^r\,,\end{equation}
(so $C^{(1)}$ is just the Carlitz module).
The associated $\tau$-sheaf is then $\underline{{\mathcal F}_g}$
for $g(\theta,T)=\frac{1}{\beta}(T-\theta)$. Thus one sees how small
a subset of all rank $1$ $\tau$-sheaves is occupied by the rank $1$
Drinfeld modules.
\end{rem}

Let $\underline{\mathcal F}$ and $\underline{\mathcal G}$ be two coherent 
$\tau$-sheaves. We define the {\it tensor product} $\tau$-sheaf $\underline{\mathcal F}\otimes
\underline{\mathcal G}$ to have underlying coherent sheaf
$\dis {\mathcal F}\otimes_{{\mathcal O}_{X\times {\mathbb A}^1}}\mathcal G$
with $\dis \tau_{_{\underline{\mathcal F}\otimes{\underline{\mathcal G}}}}:=
\tau_{_{\underline{\mathcal F}}}\otimes \tau_{_{\underline{\mathcal G}}}$.
One can check, for instance, that the tensor product of strict $\tau$-sheaves
is again a strict $\tau$-sheaf.

\begin{example}\label{carlitztensorpowers3}
Let $\underline{\mathcal C}$ be
 the $\tau$-sheaf associated to the Carlitz module, as in Example
\ref{carlitztausheaf}. One can now easily form the $n$-th tensor power
$\tau$-sheaf $\underline{\mathcal C}^{\otimes n}$. As in \cite{at1}, this
$\tau$-sheaf is isomorphic to the canonical $\tau$-sheaf associated to
$C^{\otimes n}$ which also justifies the latter's name.\end{example}

We will also identify $C^{\otimes n}$ with its associated $\tau$-sheaf
in later applications.

\begin{rem}\label{tensorproductTmodules}
Anderson \cite{an1} (also \S 5.5 of \cite{go4}) gives a
very important ``purity'' condition that insures in general that the tensor
product of the $\tau$-sheaves associated to two $T$-modules also arises
from a $T$-module.\end{rem}

\begin{rem}\label{compartTtau}
As the reader will hopefully have come to see, $T$-modules and $\tau$-sheaves
are two sides to the same coin, so to speak. 
Indeed, with $T$-modules one focuses
on the realization of $A$ as certain algebraic endomorphisms 
of ${\mathbb G}_a^e$
for some $e$. On the other hand, with $\tau$-sheaves, one emphasizes,
and generalizes, the associated motives of the $T$-modules. 
As we shall see in Subsection \ref{crystals}, it is the $\tau$-sheaf formalism that is essential in 
establishing the basic analytic properties of $L$-series in the characteristic
$p$ theory. 
However, Drinfeld modules over $\Fr((1/T))$ arise
from lattices and such lattices are needed, at least, for properties of
special-values of $L$-functions such as given in 
Example \ref{carlitzspecialvalue}. 
It is therefore natural to ask about the relationship of general
$T$-modules to lattices. In our next subsection we will discuss what
is known in this regard; it turns out that the answer is essential
for B\"ockle's theory.
\end{rem}
\subsection{Uniformization of general $T$-modules}\label{uniformT}
As before, let $k=\Fr(T)$, $\imath\colon A\to k$  the identity map and
$\theta=\imath (T)$. Let $K:=\Fr((1/\theta))$ with fixed algebraic
closure $\bar{K}$. So we are back in
the analytic set-up of Subsection \ref{drinfeld}. 
Let $E=(E_{\text gp},\psi)$ be a $t$-module of dimension $e$ defined
over a finite extension $L\subset \bar{K}$ of $K$. Without
loss of generality we can, and will, suppose that $E_{\text gp} \simeq
{\mathbb G}_a^e$. 

By definition (Equation \ref{Tmoduledef}) one knows that
$\psi_T=(\theta I_e+N)x+$\{higher terms\} with $N$ nilpotent. Clearly
the action $\psi_{T,\ast}$ of $T$ on the Lie algebra of $E$ is then
given by $\theta I_e+N$. One now formally looks for an exponential
function $\exp_E$ associated to $E$ of the form
$$\exp_E=\sum_{i=0}^\infty Q_i x^{r^i}\,$$
where $\dis x=\left(\begin{array}{c} x_1\\\vdots\\x_e\end{array}\right)\in
\text{Lie}(E)$, $x^{r^i}$ is defined in the obvious fashion, $Q_0=I_e$ and
the $Q_i$ are $e\times e$ matrices with coefficients
in $L$. As in the Drinfeld module case
(Equation \ref{basicdrinfeld}), $\exp_E$ is further required to satisfy
\begin{equation}\label{basicTmoduleexp}
\exp_E(\psi_{T,\ast} x)=\psi_T(\exp_E(x))\,.
\end{equation}
Using (\ref{basicTmoduleexp}), one readily, and uniquely, finds
the coefficient matrices $Q_i$ and that $\exp_E (x)$ is \emph{entire}
(i.e., converges for all $x$). 

However, as Anderson discovered, as soon as $e>1$ a fundamental
problem arises in that there exist abelian $T$-modules $E$ where
$\exp_E(x)$ is \emph{not} surjective on geometric points
(that is, over $\bar{K}$). Anderson \cite{an1} gives some necessary and
sufficient conditions 
for the geometric 
surjectivity of $\exp_E(x)$. We will focus here on the criterion
Anderson calls ``rigid analytic triviality.'' 

Let $M$ be the $T$-motive of $E$. Let $L\{T\}$ be the Tate algebra of
all power series $\dis \sum_{j=0}^\infty c_j T^j$ where $c_j\in L$ all $j$ and
$c_j\to 0$ as $j\to \infty$.

\begin{defn}\label{rigidm}
1.~We set $M\{T\}:=M\otimes_{L[T]}L\{T\}$ with its obvious
$L\{T\}$-module structure. We equip $M\{T\}$ with a $\tau$-action by
setting
$$\tau(m\otimes \sum c_jT^j):=\tau m\otimes \sum c_j^rT^j\,.$$
2.~We let $M\{T\}^\tau\subset M\{T\}$ be the $A$-module of $\tau$-invariants.\\
3.~The module $M$ is said to be \emph{rigid analytically trivial over $L$} if
the natural map $M\{T\}^\tau\otimes_A L\{T\}\to M\{T\}$ is an isomorphism.
\end{defn}

It is important to note that Definition \ref{rigidm} makes sense for
general arbitrary $\tau$-modules over $L$.

Anderson then proves that $\exp_E(x)$ is surjective on geometric
points if and only if there
is a finite extension $L^\prime$ of $L$ such that the motive $M$ of $E$ over
$L^\prime$ is rigid analytically trivial (over $L^\prime$). This
condition is preserved under tensor products.

If $\exp_E(x)$ is geometrically surjective, then its kernel (as a homomorphism
of groups)
$\mathcal L$ is called the ``lattice of $E$.'' One can show that
$\mathcal L$ is an $A$-module of the same rank as $M$.

\begin{defn}\label{defuniform}
We say that $E$, and $M$, is {\it uniformizable over a field $L$} if and only
if it is rigid analytically trivial over $L$. We say that $E$, and $M$,
are {\it uniformizable} if and only if there is a finite extension $L^\prime$
of $L$ over which they are uniformizable.
\end{defn}

As before, this notion can be extended to arbitrary $\tau$-modules $M$
(with no obvious exponential function attached!).

As an example, $C^{\otimes n}$ is uniformizable over $\Fr((1/\theta_1))$,
$\theta_1:=(-\theta)^{1/(r-1)}$, all $n\geq 1$, while
$C^{\otimes m(r-1)}$ is uniformizable over $\Fr((1/\theta))$ for
all $m\geq 1$. 

If $E$ is uniformizable over $L$ then the $\tau$-invariants
$M\{T\}^\tau$ form a free $A$-module of rank equal to the rank of $M$.
The converse is also true (and is an unpublished result of Urs Hartl):
If the $\tau$-invariants over $L$ form a free module of rank equal to that of
$M$, then $M$ is uniformizable over $L$. 

\begin{rem} \label{nomoreinv}
Implicit in the above statement is the assertion that if
$M\{T\}^\tau$ has rank equal to that of $M$, then one obtains
\emph{no further} invariants by passing to any finite extension $L^\prime$.
This is indeed true and can be seen directly. Indeed, the invariants
over any finite extension will have the same rank. Thus, if $m$
is one such invariant, there is an $f\in A$ such that $fm$ is an
invariant over $L$. One then sees that this forces $m$ to be defined
over $L$ also.
\end{rem}

Let $E$ be a uniformizable $T$-module which is defined
over $L$. Let $L^\prime\subset \bar{K}$ be the finite extension generated
by the lattice $\mathcal L$ of $E$. 
One then sees that $L^\prime$ is the smallest
extension of $L$ over which $E$ is uniformizable.

\begin{example}\label{ratlunifdrin}
We will present here the rank $1$ Drinfeld module $C^{(-\theta)}$ defined
over $\Fr(\theta)$ by
\begin{equation}\label{ctheta}
C^{(-\theta)}_T(x)=\theta x-\theta x^r\,. \end{equation}
Using the explicit knowledge of the period $\xi$ of the Carlitz module
(see Example \ref{carlitz1}), one sees readily that the lattice of
$C^{(-\theta)}$ lies in $\Fr((1/\theta))$; thus $C^{(-\theta)}$ is
uniformizable over $\Fr((1/\theta))$.
\end{example}

\subsection{Tate modules of Drinfeld modules and $T$-modules
over $\Fr(\theta)$}\label{tatemodules}
One approach to constructing the $L$-series of an elliptic
curve over $\Q$ mentioned in Subsection \ref{lseriesellcurve}
is through the use of its Tate-modules. 
We will use the same approach here to define the $L$-series of
general Drinfeld modules and $T$-modules and, in the next subsection,
we will present the construction for $\tau$-sheaves.

Thus let $E=(E_\text{gp},\psi_E)$ be an abelian $T$-module defined over $k=\Fr(\theta)$.
Let $v=(g)$ be the prime associated to a monic irreducible $g\in A$. We define
the $v^i$ torsion points of $E$ to be the kernel of the map
$x\mapsto \psi_{E,g^i}(x)$ where $x\in E_\text{gp}(\bar{k})$ and
$\bar k$ is a fixed algebraic closure of $k$; we denote
this kernel by ``$E[v^i]$.'' Clearly, $E[v^i]$ inherits an $A$-structure
and it can then be shown that $E[v^i]\simeq (A/v^i)^t$ where $t$ is the
rank of $E$. The $v$-adic Tate module of $E$ is then defined by
\begin{equation}\label{Tatemodule}
T_v(E):=\varprojlim_i E[v^i]\,.
\end{equation}
Thus $T_v(E)$ is a free $A_v$-module of rank $t$.
Finally, we set 
\begin{equation}\label{H1Tmodule}
H^1_v(E)=\text{Hom}_{A_v}(T_v(E),k_v)\,.\end{equation}

The various $A_v$-modules, $\{H^1_v(E)\}$, form a compatible system of 
Galois
representations as with elliptic curves. Using geometric Frobenius elements and
invariants of inertia, again as in the case of elliptic curves, one
obtains local $L$-factors $L_f(E,u)$ for monic primes $f\in A$ with
$L_f(E,u)^{-1}\in A[u]$. One then defines the $L$-function of $E$,
$L(E,s)$ for $s\in S_\infty$, by
\begin{equation}\label{LTmodule}
L(E,s):=\prod_f L_f(E,f^{-s})\,.
\end{equation}
One sees easily that $L(C,s)$, with the above definition, agrees with
$L(C,s)$ as given in (\ref{Lcarlitz1}).

Two $T$-modules are said to be isogenous if there is a finite surjective
map between them (i.e., a surjective map of the underlying algebraic
groups which commutes with the $A$-actions). It is known that the
isogeny class for Drinfeld
modules and many $T$-modules (\cite{tag1}, \cite{tag2},
\cite{tam1}) is 
determined by the associated $L$-series (as one can read off the local
$L$-factors from the $L$-series).

In \cite{ga2} Gardeyn shows that an abelian $T$-module is uniformizable
if and only if the Tate action of the decomposition group at $\infty$ has 
finite image.
\subsection{The $L$-series of a $\tau$-sheaf over $k$}\label{Ltau}
As is discussed in \cite{ga1}, general $T$-modules are not the proper setting
in which to analyze the local factors of their associated $L$-series. It is
relatively easy to define the appropriate notion of ``good'' prime
for a $T$-module $E$ (one just wants to insure that one can reduce
the $T$-action of $E$ to obtain
a $T$-module of the {\it same} rank over the quotient field.) 
However, outside of the case of
Drinfeld modules, one then loses the connection between good primes
for $E$ and good (=unramified) primes for the compatible system
$\{H^1_v(E)\}$. Moreover, even in the case of Drinfeld modules,
$\tau$-sheaves are needed in order to describe the factors at the
bad primes, see Example \ref{francis3}.

The techniques for defining the $L$-function of a $\tau$-sheaf goes back to
work of Anderson \cite{an1} 
on $T$-modules. Let $E$ be an abelian $T$-module with
associated motive $M=M(E)$ as in Subsection \ref{taumodule} and let
$\bar{M}$ be constructed in the same fashion as $M$ but over the algebraic
closure $\bar k$ of $k$.
Let $v=(g)$ be a prime of $A$. Then
Anderson shows that the Galois module $H^1_v(E)$ is isomorphic to the
Galois module $H^1_v(M)$ defined by
\begin{equation}\label{H1motive}
H^1_v(M):=\varprojlim_i (\bar{M}/g^i\bar{M})^\tau\,,\end{equation}
(where $N^\tau:=\{\lambda \in N \mid \tau\lambda =\lambda\}$ for any $\tau$-module $N$).

The above definitions immediately carry over to the case of $\tau$-sheaves
$\underline{\mathcal F}=({\mathcal F},\tau)$
over $k$; one obtains local factors $L_f(\underline{\mathcal F},u)$ again
using inertial invariants and characteristic polynomials of geometric
Frobenius elements.
The idea of G.\ B\"ockle, R.\ Pink and F.\ Gardeyn (again
following work of Anderson), is to
show that $L_f(\underline{\mathcal F},u)$
be expressed in terms of the $\tau$-action
itself. Indeed, at a bad prime $f$ (where
there are non-trivial invariants of inertia) Gardeyn \cite{ga3} constructs a
``maximal model'' $\underline{\mathcal F}^M=({\mathcal F}^M,\tau^M)$ 
of $\underline{\mathcal F}$ (which may be viewed as a ``N\'eron model''
for $\underline{\mathcal F}$). The point is that the special fiber
$\underline{\mathcal F}_\text{sp}=({\mathcal F}_\text{sp},\tau_\text{sp})$ of 
$\underline{\mathcal F}^M$ is a $\tau$-sheaf on $\text{Spec}(\F_f)$ (where
$\F_f$ is the residue field at $f$). One then sees that 
\begin{equation}\label{gebhard}
L_f(\underline{\mathcal F},u)^{-1}=\det_A\left(1-u\tau\mid H^0({\mathcal F}_\text{sp})\right)\,,\end{equation}
which establishes, for instance, that 
$L_f(\underline{\mathcal F},u)^{-1}\in A[u]$.

\begin{example}\label{francis3}
Let $\psi$ be a Drinfeld module over $k$. In \cite{ga1}, Gardeyn presents the
local factors $L_f(\psi,u)$ at the bad primes $f$. It is shown that if
$\psi$ has bad reduction at $f$ but $\psi$ obtains good reduction over a finite
extension $L$ of $k$ (and a prime of $L$ above $f$) 
then $L_f(\psi,u)=1$. Moreover, if there does not exist
a finite extension $L$ of $k$ over which $\psi$ obtains good reduction, then
$L_f(\psi,u)^{-1}\in \Fr[u]\subset A[u]$. This is remarkably similar to the
case of elliptic curves (\ref{ellipticfactor2}). It would be interesting to
establish exactly which polynomials in $\Fr[u]$ actually occur for
a given Drinfeld module $\psi$. Note also that all rank $1$ Drinfeld modules
have potentially good reduction (since they are all geometrically isomorphic
to the Carlitz module). As such, the local factors at the bad primes in the
rank $1$ case are all identically $1$ as one expects.\end{example}

\begin{rem}\label{gebhard2}
In \cite{boc1}, the local $L$-factors of a $\tau$-sheaf $\underline{\mathcal F}$ are
defined directly as in Equation (\ref{gebhard}) without using Galois 
representations. However, to any $\tau$-sheaf one can attach a
constructible \'etale sheaf of $A_v$-modules which is a natural Galois
module. One can use this Galois module as we have used $H^1_v(E)$ for
a $T$-module to define the $L$-factor (and, indeed, in the $T$-module case
the Galois module is isomorphic to $H^1(E)$). Therefore one can always
use the classical Galois formalism to define local $L$-factors in general.
\end{rem}

We will finish this subsection by describing briefly the Galois representations
associated to $\tau$-sheaves $\underline{{\mathcal F}_g}=(\mathcal{F}_g,\tau_g)$
where $\mathcal{F}_g=\text{Spec}(\Fr(\theta)[T])$, 
$0\neq g(\theta,T)\in \Fr(\theta)[T]$, and $\tau_g$ is given by
(\ref{carlitztausheaf3}). As these sheaves have rank $1$, we obtain
$1$-dimensional $v$-adic representations which we
denote by $\rho_{g,v}$. Let $f(\theta)\in \Fr[\theta]$ be a monic irreducible
polynomial of degree $d$ with roots $\{\bar{\theta},\bar{\theta}^r,
\ldots, \bar{\theta}^{r^{d-1}}\}$. Set
\begin{equation}\label{gebhard3}
g^f(T):=\prod_{i=0}^{d-1} g(T,\bar{\theta}^{r^i})\in \Fr[T]\,.
\end{equation}
For instance if $g(T,\theta)=\prod_i (h_i(T)-\theta)$, where $h_i(T)$ does not
involve $\theta$, then $g^f(T)=\prod_if(h_i(T))$.
Suppose now that $g^f(T)\in A_v^\ast$ and, finally, let
$\text{Frob}_f$ be the geometric Frobenius at $(f(\theta))$. Then one
has 
\begin{equation}\label{boecklegalois}
\rho_{g,v}(\text{Frob}_f)=g^f(T)\,.\end{equation}
(I am indebted to B\"ockle for
pointing out this simple and elegant formula.)

Now let $C$ be the Carlitz module. One knows
that $C$ corresponds to the function $g(\theta,T)=T-\theta$.
Let $v$ be as above and denote $\rho_{g,v}$ by $\rho_{C,v}$.
Let $f(\theta)$ be a monic prime with $v$ relatively prime to
$f(T)$. Then one has 
\begin{equation}\label{boecklegalois2}
\rho_{C,v}(\text{Frob}_f)=f(T)\in A_v^\ast
\end{equation}
which agrees with (\ref{charpL1}) and where we recall that we use the dual
action to define $L$-series.

More generally let $C^{(\beta)}$ be the general rank $1$ Drinfeld
module over $\Fr(\theta)$ as in Remark \ref{taumodulerankone} with
associated function $g(\theta,T)=\frac{1}{\beta}(T-\theta)$. Let
$\rho_{C^{(\beta)},v}$ be the associated $v$-adic representation.
Then, as $\beta$ is \emph{constant} in $T$ (so that $\beta^f(T)$ is also
constant in $T$), one finds
\begin{equation}\label{boecklegalois3}
\rho_{C^{(\beta)},v}=\chi_\beta \rho_{C,v}\end{equation}
where $\chi_\beta$ is an $\Fr^\ast$-valued Galois
character which is independent of $v$.

\subsection{Special polynomials and Carlitz tensor powers}\label{specialcarlitz}
Recall that in the case of the $L$-series $L(C,s)$ of the Carlitz module the
functions $L(C,x/T^i,-i)$, $i$ a non-negative integer, actually
belong to $A[x^{-1}]$. Let $\underline{\mathcal F}$ now be a $\tau$-sheaf
with $L$-series $L(\underline{\mathcal F},s)$. The case of the Carlitz
module suggests looking at the power series $L(\underline{\mathcal F},x/T^i,-i)$
for $i$ as above. In our next subsection we will establish that these
{\it special power series} are in fact rational functions
 with $A$-coefficients (and,
naturally, called the {\it special functions of 
$L(\underline{\mathcal F},s)$}). Essential
to the proof is the equality
\begin{equation}\label{tensorwithcarlitz}
L(\underline{\mathcal F},x/T^i,-i)=
L(\underline{\mathcal F}\otimes C^{\otimes i},
x,0)\,.\end{equation}
Equation (\ref{tensorwithcarlitz}) follows directly from looking at
the associated Galois representations. In particular, one concludes for non-negative
integers $i$ that
\begin{equation}\label{tensorwithcarlitz2}
L(\underline{\mathcal F},s-i)=L(\underline{\mathcal F}, s-s_i)=
L(\underline{\mathcal F}\otimes C^{\otimes i},s)\,.\end{equation}

\subsection{Crystals and their cohomology}\label{crystals}
In this subsection we review briefly the theory of ``crystals'' associated
to $\tau$-sheaves developed by R.\ Pink and G.\ B\"ockle 
\cite{bp1} (see also \cite{boc1} and \cite{boc2}). Let
$\underline{\mathcal F}=(\mathcal{F},\tau)$ be a $\tau$-sheaf on
a scheme $X$ where $\tau$ acts nilpotently, that is, $\tau^m=0$ for
some $m>0$. From Equation (\ref{gebhard}) we see that the 
$L$-factors of $\underline{\mathcal F}$ will be trivial 
(identically $1$) at every prime;
thus the associated global $L$-series will also be trivial. Therefore,
from the $L$-series point of view, such $\tau$-sheaves are negligible. The
idea of B\"ockle and Pink is to make this precise by passing to a certain
quotient category. 
More precisely, the category of ``crystals over $X$'' 
is the quotient category of $\tau$-sheaves modulo the subcategory of 
nilpotent $\tau$-sheaves.

The category of crystals has the advantage that a cohomology theory
may be developed for it. This cohomology theory is very closely related to
coherent sheaf cohomology but which possesses only the first three of
the canonical six functors $\{Rf_!, f^\ast,\otimes,f_\ast, f^!,\text{Hom}\}$.
However, the cohomology of crystals does possess a Lefschetz trace formula.
As such, by using (\ref{tensorwithcarlitz}), B\"ockle and Pink establish
that the special power series associated to $L(\underline{\mathcal F},s)$
are rational functions. If, for instance, $\underline
{\mathcal F}$ is locally free, then one obtains an entire function
whose special rational functions are polynomials (Th. 4.15 of \cite{boc1}).
Furthermore,  B\"ockle establishes that the
degrees of these special polynomials $L(\underline{\mathcal F},x/T^i,-i)$
grow \emph{logarithmically} in $i$. This is
then enough to establish that general $L(\underline{\mathcal F},s)$  
have meromorphic interpolations at all the places of $k$. 

\subsection{Modular forms in characteristic $p$}\label{charpmodforms}
We now have all the techniques necessary to begin studying modular forms
in characteristic $p$ which we present in this subsection. Let
$\Omega$ be the Drinfeld upper-half plane as given in Definition
\ref{omega2}. Based on the discussion given in Subsection \ref{modularforms},
the notion of a ``congruence subgroup'' $\tilde{\Gamma}$ of
$\Gamma:=GL_2(A)$ is obvious as is the notion of an unrestricted
modular form of weight $j$ (where $j$ is now an integer) for $\tilde{\Gamma}$ 
(upon replacing ``analytic'' with ``rigid analytic'' in Definition
\ref{unrestmod}). 

Thus, following classical theory, we clearly
need to describe what happens at
the cusps $\tilde{\Gamma}\backslash {\mathbb P}^1(k)$ and to do this one 
needs only
treat the special case of the cusp $\infty$. As before let 
$\tilde{\Gamma}_\infty$
be the subgroup of $\tilde{\Gamma}$ that fixes $\infty$. One sees that 
$\tilde{\Gamma}$ consists of mappings of the form $z\mapsto \alpha z+b$ where
$\alpha$ belongs to a subgroup $H$ of $\Fr^\ast$ and $b\in I$ where $I=(i)$ is 
an ideal of $A$. We set $e_\infty(z):=e_C(\xi z/i)$ where $e_C(z)$ is the 
exponential
of the Carlitz module and $\xi$ is its period. Finally we set
$q:=e_\infty(z)^{-e}$ where $e$ is the order of $H$. In \cite{go1}, 
it is shown that
$q$ is a parameter at the cusp $\infty$.

With the above choice of parameter $q$, the definitions of 
{\it holomorphic form} and
{\it cusp form} are exactly the same as their complex counterparts. One can
show (ibid.) that holomorphic forms are sections of line bundles on the
 associated
compactified moduli curves; therefore for fixed subgroup and weight, they form
finite dimensional $\bar K$-vector spaces.

There are now $2$ distinct cases of subgroups $\tilde{\Gamma}$ of
interest to us.
The first case is $\tilde{\Gamma}=\Gamma$ and the
second case is the full congruence subgroup
$$\tilde{\Gamma}=\Gamma(N)=\left\{\gamma\in GL_2(A)\mid 
\gamma \equiv \begin{pmatrix} 1&0\\0&1\end{pmatrix}\pmod{N}\right\}$$ 
for some polynomial $N\in A$. In the first case, the parameter $q$
at $\infty$ is
$e_C(\xi z)^{1-r}$ and in the second case it is $q=e_C(\xi z/N)^{-1}$.

\begin{rem}\label{dq}
In the classical elliptic modular theory one has $dq=cq\cdot d\tau$ 
for some non-zero
constant $c$. Thus one sees that cusp forms of weight $2$ correspond to 
holomorphic differential forms on the associated complete moduli curve. 
For $\Gamma(N)$, with $N\in A$, one computes readily that $dq=cq^2\cdot dz$
with $c\neq 0$. Thus cusp forms
with zeroes of order $2$ at every cusps correspond to holomorphic differential
forms on the associated complete modular curve. Such cusp forms are
called ``double cusp forms.'' 
\end{rem}

We denote the space of cusp forms of weight $j$ associated to $N$ by $S(N,j)$ and the
subspace of double cusp forms by $S^2(N,j)$. A simple calculation implies that
a cusp form $f$ for $GL_2(A)$ {\it automatically} 
becomes a double cusp form for
$\Gamma(N)$ for any non-constant $N$.

\begin{rem} \label{type}
 Recall that after Definition \ref{unrestmod} we
mentioned ``multiplier systems'' that allow one to obtain a (slightly)
generalized notion of modular forms. One such
multiplier is $\det \left(\begin{array}{cc} a&b\\c&d\end{array}\right)^{-t}=
(ad-bc)^{-t}$ where $t$ is an integer; one then says that the modular form 
has {\it type} $t$ (see, e.g., \cite{ge1} or Definition 5.1 of \cite{boc2}).
\end{rem}

\begin{example}\label{ganddelta} 
As in (\ref{rank2drinfeld}), a rank two Drinfeld $A$-module $\psi$ is
uniquely determined
by $\psi_T(x)=Tx+g(\psi)x^r+\Delta(\psi) x^{r^2}$ where $\Delta(\psi)\neq 0$. 
Let $z\in\Omega$
and let $L_z:=A+Az$ be the associated rank $2$ $A$-lattice and
$\psi^{(z)}$ the associated rank $2$ Drinfeld module. As in the classical
case, the maps
$g\colon z\mapsto g(\psi^{(z)})$ and $\Delta\colon z\mapsto \Delta(\psi^{(z)})$
define rigid analytic modular forms for the group $GL_2(A)$ of
weights $r-1$ and $r^2-1$ respectively (and type $0$). 
Moreover, $\Delta$ is easily
seen to be a cusp form as it is classically.\end{example} 

\begin{rem}\label{caveat1}
When working with modular forms associated to congruence subgroups
there is a major difference between classical theory and the
theory developed in \cite{boc2}. Indeed, B\"ockle needs to work with
the {\it full} moduli spaces attached to congruence subgroups and in
particular the moduli space of Drinfeld modules of rank $2$ with
level $I$ structure (i.e., a basis for the $i$-division points). This
moduli space contains {\it many} different geometric components (as does
its classical counterpart). As in Drinfeld's original paper \cite{dr1},
these components are best handled through the use of the adeles. In particular,
in \cite{boc2}, \S 5.5, B\"ockle develops a theory of types for
adelic modular forms which generalizes that given in Remark \ref{type}
above.\end{rem}

The definition of the Hecke operators $T(I)$ for ideals
$I$ of A is then modeled on
classical theory. In particular B\"ockle \cite{boc2} presents
naturally defined Hecke-operators in the adelic setting (and for
general base rings $A$) which depend on the type, weight and
level involved. Moreover, these Hecke operators do {\it not}
fix the components of the underlying moduli spaces. In particular they
therefore differ from the ones defined in \cite{go1}, \cite{go2}
and \cite{ge1}; the latter Hecke operators fix the components of
the moduli spaces but cannot be defined for general base rings
$A$. When recalling the Eichler-Shimura isomorphism given in
\cite{boc2}, this is an important consideration; a comparison
between the two viewpoints is given in Example 6.13 of \cite{boc2}.

One sees naturally that the cusp forms of a given weight
are stable under the Hecke operators, but also, when $\tilde{\Gamma}=
\Gamma (N)$, so are the double cusp forms (which is highly remarkable
from the classical viewpoint!). Moreover, as in the classical case,
the Hecke operators form a commutative ring of endomorphisms of these
spaces.
 
Thus there are really three Hecke stable spaces of interest: $S(N,j)$,
$S^2(N,j)$ and the quotient space $S(N,j)/S^2(N,j)$.

\begin{rem}\label{heckeisstrange}
Classical Hecke operators, as covered in Subsection \ref{modularforms}, have the
property that $T(p^2)\neq T(p)^2$ for $p$ prime. Remarkably, in the
characteristic $p$ case, one finds that $T(I^2)=T(I)^2$ for {\it any}
ideal $I$ (including, precisely, the case of $I$ prime);
thus the Hecke operators are {\it strongly multiplicative}. Indeed, the
classical computation of $T(\wp^2)$, $\wp$ prime, works and one sees that the
terms different from $T(\wp)^2$ are weighted with integer factors divisible
by $p=0\in \Fr$. This commutativity is essentially
the reason that the Hecke operators give rise to {\it abelian} representations
as in our next subsection.
\end{rem}

\subsection{Galois representations associated to cusp forms}\label{galrepmodforms}
In this subsection we summarize very briefly the results of \cite{boc2}
on Galois representations associated to cusp forms in characteristic $p$.
 Fix $N\in A$ and view $S(N,j)$ as a Hecke-module. As the ring of Hecke
operators is commutative we can decompose $S(N,j)$ into generalized
eigenspaces. Let $\{M_1,\ldots,M_\lambda\}$ denote the simple
Hecke subfactors of the {\it true} eigenspaces. Every simple Hecke
subfactor of $S(N,j)$ is then isomorphic to one of the $M_i$.

To each $M_i$ corresponds to a true Hecke eigenform $f_i\in S(N,j)$. Let
$\mathfrak P$ be a prime of $A$ not dividing $N$ and suppose that 
$T({\mathfrak P})f_i=\alpha_{i,\mathfrak P} f_i$ (where
we use the adelic Hecke operators of \cite{boc2}). Via the general 
cohomological formalism of
crystals, B\"ockle attaches to each $f_i$ 
a rank $1$ $\tau$-sheaf $\underline{\mathcal M}_i$; this is done in a 
non-canonical fashion. 

Let $v$ be a prime of $A$. The general theory of $\tau$-sheaves, as in
Remark \ref{gebhard2},
then gives us a continuous $1$-dimensional $v$-adic Galois
 representation $\rho_i=\rho_{f_i}$
for each $i$ (which {\it is} indeed canonical!). We call the compatible
system of representations obtained this way the {\it B\"ockle system 
(of Galois representations) attached to $f_i$}.
The Eichler-Shimura relation established in \cite{boc2}
in this context then implies that
$$\rho_i(\text{Frob}_{\mathfrak P})=\alpha_{i,\mathfrak P}$$
for $\mathfrak P$ prime to $N$ and $v$. In particular, we conclude that
$\alpha_{i,\mathfrak P}\neq 0$.

\begin{rem}\label{bocklewarning}
In \cite{boc2}, Theorem 13.2, the above result is only established for
cusp forms of weight $n$, type $n-1$ and level $I\neq A$. There is also
a general ``yoga'' which allows one to change types, after increasing
the level; by using compatibilities of the Galois representations
attached to modular forms, the above result can be extended to
arbitrary types independent of the weight, cf.\ \cite{boc2}
Lemma 5.32 and Remark 6.12. This is important to us since we want to
consider cusp forms of level $A$ (attached the full modular group)
and type $0$. A more conceptual proof which avoids this yoga
can be given by proving an Eichler-Shimura isomorphism
for fixed level and arbitrary weight $n$ for all types $l\geq n-1$,
where however the $\tau$-sheaf $\underline{\mathcal M}_i$ mentioned above have
to be twisted suitably. ``Untwisting'' by powers of $C^{\otimes (r-1)}$
on the Galois side, one then obtains the result.

While the above process attaches Galois representations to cusp forms
for any type $l$, for $l<n-1$ there is {\it no} $\tau$-sheaf
associated to the representation. This is similar to the classical
situation where the inverse of the Tate motive is not represented
by a geometric object but the corresponding cyclotomic character
obviously has an inverse.
\end{rem}

\begin{rem}\label{allfinite}
If the cusp form is {\it not} double-cuspidal, then the associated
family of Galois representations arises essentially
from a finite character. To be more precise, if the weight
is $n$ and the type is $n-1$, then there is indeed a finite character.
For other types and the same weight, the Galois representations
get twisted by some natural $1$-dimensional characters associated to
Drinfeld modular varieties of rank $1$ Drinfeld modules.
Moreover,
B\"ockle establishes that the class of finite characters which
arise are \emph{all} finite characteristic $p$
valued characters allowed by the explicit class field theory
of rank $1$ Drinfeld modules. That is one obtains those
finite dimensional characters of abelian extensions of $k=\Fr(\theta)$
which are \emph{totally-split} at $\infty$. Moreover,
one obtains such representations for arbitrary
weights $>2$ (for $2$ there are some modifications involving
the trivial character). It is reasonable to view
the associated cusp forms as being rather analogous to the theta-series
assigned to finite characters classically as in Equation (\ref{theta3}).
Explicitly constructing such cusp forms in characteristic $p$
is certainly now a very interesting problem.
\end{rem}

\begin{rem}\label{Lfiisentire}
Although the choice of $\tau$-sheaf 
$\underline{\mathcal M}_i$ associated to $f_i$ is not
canonical, the associated Galois representations are and depend only
on the Hecke eigenvalues. As such, one can define 
$L(f_i,s):=L(\underline{\mathcal M}_i,s)$ in an unambiguous
fashion. The results of B\"ockle in Subsection \ref{crystals} then imply
the analytic continuation of $L(f_i,s)$ (at all places of $k$).
\end{rem}

\begin{question}\label{twistchar}
Classically one can twist modular forms by characters simply by multiplication
of the Fourier coefficients. Can one define such twists for the
B\"ockle systems in characteristic $p$?
\end{question}

\begin{rem}\label{needuniform}
It is very important for us that
B\"ockle's theory \emph{does} establish at least
one (so far!)
constraint on the $\tau$-sheaves that may arise from modular forms.
Indeed, B\"ockle shows that such $\tau$-sheaves arise from decomposing
a $\tau$-sheaf (via ``complex multiplications'') which is 
defined and \emph{uniformizable} (in the sense of Definition
\ref{defuniform}) over $\Fr((1/\theta))$.
\end{rem}

\begin{rem}\label{whygamma1} The reader may well be asking why
one works with $\Gamma(N)$ as opposed to $\Gamma_0(N)$.
One does not use $\Gamma_0(N)$
because one needs a fine moduli space (i.e., a representable functor)
for B\"ockle's constructions
and $\Gamma_0(N)$ is not associated with a representable functor.
Indeed, B\"ockle begins with the $\tau$-sheaf $\underline{\mathcal F}_N$
 associated to the
universal family of Drinfeld modules associated to $\Gamma(N)$. The 
representations arise by relating the 
B\"ockle-Pink cohomology of the symmetric powers of $\underline{\mathcal F}_N$
with \'etale cohomology.
\end{rem}

\subsection{Modularity for rank $1$ Drinfeld modules}\label{rank1taumod}
Let $\underline{\mathcal F}$ be any $\tau$-sheaf defined over $k$ and let
$L(\underline{\mathcal F},s)$, $s\in S_\infty$, be its $L$-series. From
Equation (\ref{tensorwithcarlitz2}) we see that the $L$-series of
$\underline{\mathcal F}$ and $\underline{\mathcal F}\otimes C^{\otimes n}$
are simple integral {\it translates} of each other.
 Thus, from the point of view
of $L$-series, the $\tau$-sheaves 
$\underline{\mathcal F}$ and $\underline{\mathcal F}
\otimes C^{\otimes n}$ are equivalent.

The above observation will guide our definition of ``modularity.''
In fact, there are really two notions of ``modularity'' implicit
in the theory. Let $k=\Fr(\theta)$ as before and let $k^\text{sep}$
be a fixed separable closure.

\begin{defn}\label{modularity}
We say that a Drinfeld module $\psi$ defined over $k=\Fr(\theta)$
is {\it modular of class I} if and only if its $L$-series is an
integral translate of the
$L$-series of a finite $\Fr^\ast$-valued character of
$\text{Gal}(k^\text{sep}/k)$ which has trivial component
at $\infty$ (see Remark \ref{allfinite}). We say $\psi$ is {\it modular of 
class II} if and only if its $L$-series is an
integral translate of $L(f,s)$ where
$f$ is a double cusp form of some weight and level.
\end{defn}

\begin{rem}\label{modularrankone}
In general for a Drinfeld module $\phi$ to possibly have its
Galois representations (possibly twisted by those of $C^{j(r-1)}$)
arise from one of the $\tau$-sheaves ${\mathcal M}_i$, it must
have {\it abelian} Galois representations on its Tate modules. Thus it
is either of rank $1$ or has ``complex multiplication'' when the rank
$d>1$. Since the Galois-image is abelian, the latter means that the ring
$A^\prime$ of endomorphism of $\phi$ is commutative and a finite
extension of $A$. Furthermore by Prop.\ 4.7.17 of \cite{go4}
$k^\prime:=A^\prime \otimes_A k$ as well as $K^\prime:=k^\prime \otimes_k
K$ (where we recall $K=k_\infty=\Fr((1/T))$\,) are fields. Since $\phi$ must
be uniformizable over $K$, it therefore must also satisfy the
weaker condition that it is uniformizable over $K^\prime$.
\end{rem}

\begin{example}\label{typeIexamples}
Let $C^{(\beta)}$ be the general rank $1$ Drinfeld module given
in Equation (\ref{taumodulerankone2}) with $\beta\in\Fr(\theta)$. 
Suppose that $\beta=\alpha^{r-1}$ with $\alpha \in \Fr((1/\theta))$;
thus over $\Fr((1/\theta))$ one has $C^{(\beta)}_a(x)=\alpha^{-1}x\circ
C_a(x) \circ \alpha x$. In particular, $C$ and $C^{\beta}$ are isomorphic
over $\Fr((1/\theta))$. Thus the character $\chi_\beta$ of
Equation (\ref{boecklegalois3}) has trivial component at $\infty$ and
$C^{(\beta)}$ is modular of class I.
As an example, one can take $\beta:=\frac{\theta+1}{\theta}=1+
\frac{1}{\theta}$
and then find $\alpha$ via the binomial theorem applied to $1/(r-1)$.
\end{example}

\noindent
In particular, the Carlitz module is obviously then modular of class
I. We now show how it is also modular of class II.

\begin{example}\label{carlitzismodular}
Let $\Delta$ be as in Example \ref{ganddelta}.
It is shown in \cite{go2} that if ${\mathfrak P}=({\mathfrak p})$ then
\begin{equation}\label{eigenvaluesdelta}
T({\mathfrak P})\Delta={\mathfrak p}^{r-1}\Delta\,,\end{equation}
where we have used the Hecke operators from \cite{go2}.
If instead we had used the Hecke operators as in \cite{boc2}, 
we would obtain
\begin{equation}\label{eigenvaluesdeltaalboeckle}
T({\mathfrak P})\Delta={\mathfrak p}^{(r-r^2)}\Delta\,,
\end{equation}
cf.\ Example 6.13 of \cite{boc2}.
Thus the $L$-function of the
B\"ockle system of Galois representations associated to $\Delta$
equals $L(C, s+r^2-r+1)$. In particular
$C$ is therefore modular of class II. \end{example}

It should be pointed out that, in line with the results mentioned in
Subsection \ref{uniformT}, $C^{\otimes (r-1)}$ is actually uniformizable
over $\Fr((1/\theta))$.

Note that in particular, $\zeta_A(s+1-r)=L(\Delta,s)$. This should be compared
with the classical formula of Subsection \ref{thetadir} where the
theta function $\theta(\tau)$ naturally gives $\zeta(2s)$
(i.e., one needs the factor $s/2$ in the integral (\ref{zeta1})). 

In fact, in \cite{go2} the exact {\it same} result (\ref{eigenvaluesdelta})
is also established for the cusp form $g^r\Delta$ of 
weight $(2r+1)(r-1)$ (where $g$ is also defined in 
Example \ref{ganddelta}).  From
the classical viewpoint this is {\it highly} surprising!

Recall that in Example \ref{ratlunifdrin} we discussed the rank one
Drinfeld module $C^{(-\theta)}$ which is uniformizable over
$\Fr((1/\theta))$. There is \emph{no} known obstruction for the 
$v$-adic representations associated to $C^{(-\theta)}$ 
(or any Drinfeld module defined over $\Fr(\theta)$ which is uniformizable
over $\Fr((1/\theta))$\,) to be the B\"ockle system arising from some
double cusp form. We are thus led to the following question.

\begin{question}\label{q1}
Does there exist a double cusp form of some weight and level whose
B\"ockle system of Galois representations is the same as the
system arising from $C^{(-\theta)}\otimes C^{\otimes j}$ for some $j\geq 0$ 
with $j\equiv 0\pmod{r-1}$?\end{question}

\noindent
In other words, is $C^{(-\theta)}$ modular of class II with
the integer giving the translation being divisible by $r-1$?

\begin{rem}\label{refinement}
Recall that we defined the finite Galois character $\chi_\beta$
associated to $C^{(\beta)}$ in Equation (\ref{boecklegalois3}). 
When $\beta=-\theta$ standard calculations involved in the
$(r-1)$-st power reciprocity law for $\Fr[T]$ tell us that
the finite part of the
conductor of $\chi_{-\theta}$ (in the usual sense of class
field theory) is $(T)$ (see, e.g., the Exercises to \S 12 of
\cite{ros1}). Thus the level in Question \ref{q1} should almost
certainly be $(T)$. Predicting the weight is much more difficult
as the functor $Rf_\ast$ on crystals does not preserve purity;
thus there is as yet no obvious guess for the weight.
\end{rem}

The answers to Question \ref{q1} and its
refinement (Remark \ref{refinement}), as well as the analogous questions for
arbitrary rank $1$ Drinfeld modules over $\Fr(\theta)$ uniformizable over
$\Fr((1/\theta))$, will be very interesting. Classical theory leads us
to expect ``good'' reasons for the answer whether affirmative or negative.

\subsection{Final remarks}\label{finrem}
There are any number of interesting problems and comments
that virtually leap at one
from B\"ockle's constructions. We mention just a few here.

The first obvious problem is to characterize the ``Dirichlet series'' that
arise from cusp or double-cusp forms; i.e., what special properties does 
$L(f,s)$ possess besides an analytic continuation (which, after all, exists
for {\it all} $L$-series of $\tau$-sheaves)? Classically, such information
is contained in the functional equations satisfied by the Dirichlet series.
Moreover there are some questions \cite{go6} about the zeroes
of the characteristic $p$ functions that seem to be
quite natural. Furthermore, the analogy with classical theory would suggest that
the answer to these questions would involve some sort of ``functional equation'' 
in the characteristic $p$ theory. However, at present, one does not know 
even how to guess at the formulation of such a functional
equation. 
  
Secondly, the example of $\Delta$ and $g^r\Delta$ as well as the
results of B\"ockle mentioned in Remark \ref{allfinite} shows that the
relationship between the $q$-expansion of an eigenform and its 
Hecke eigenvalues is \emph{very} different
from that known classically.
In fact, one does not yet have
formulae which allow one to characterize the $q$-expansion coefficients
from the Hecke eigenvalues. 
An obvious problem is to
find additional structure that allows one to distinguish the different
cusp forms which have the same $L$-function (or, even,
the same up to translation).  As of now there is no obvious guess here also.

One would also like an explicit basis of eigenforms for the complement of
the space of double cusp forms in the space of cusp forms. There are
examples in \cite{boc2} where such forms are given by Poincar\'e series,
but no general construction is now known.

The theory of Drinfeld modules exists in the very general set-up where
$A$ can be the ring of functions in any global field $k$
of finite characteristic
regular away from a fixed place $\infty$. Virtually all of the theory
discussed above goes over directly in this general set-up. However,
when $A$ is not factorial, the reader should keep in mind that
there are NO Drinfeld modules defined over $k$ itself; rather one must
work over some Hilbert class field.

Finally the theory of rigid modular forms exists for Drinfeld modules
of {\it all} ranks. In the case $A=\Fr[T]$ there is a compactification
for these general moduli schemes of arbitrary rank
due to M.\ Kapranov \cite{ka1}. In \cite{go3}, it is
shown that Kapranov's compactification, and coherent
cohomology, allow one to conclude the
finite dimensionality of spaces of modular forms in general. It is
very reasonable to expect that B\"ockle's techniques will also work
here too, thus producing another huge class of rank $1$ $\tau$-sheaves
which will also need to be understood and somehow classified.

\end{document}